\documentclass[1p,times,review]{elsarticle}

\usepackage{hyperref}

\usepackage{numcompress}

\usepackage{amsmath,amsthm}
\usepackage{amssymb}
\usepackage{graphicx}
\usepackage{psfrag}
\usepackage{float}
\usepackage{epstopdf}
\usepackage{graphics}
\usepackage{epsf}
\usepackage{xfrac}
\usepackage{amscd}
\usepackage{color} 
\usepackage{wrapfig}
\usepackage[mathscr]{euscript}
\usepackage[english]{babel}
\usepackage{epsfig}
\usepackage{dsfont}
\usepackage[center]{caption}
\everymath{\displaystyle}
\usepackage{multirow}
\usepackage{t1enc}
\usepackage[utf8]{inputenc}
\usepackage{amsfonts}
\usepackage{stmaryrd}

\graphicspath{{AdvDiffFig/}}

\definecolor{freeblue}{rgb}{0.25,0.41,0.88}
\definecolor{darkorange}{rgb}{1.0, 0.55, 0.0}
\definecolor{mediumred-violet}{rgb}{0.73, 0.2, 0.52}
\definecolor{americanrose}{rgb}{1.0,0.01,0.24}
\definecolor{deepmagenta}{rgb}{0.8, 0.0, 0.8}
\definecolor{palatinateblue}{rgb}{0.15, 0.23, 0.89}
\definecolor{ruby}{rgb}{0.88, 0.07, 0.37}
\definecolor{shamrockgreen}{rgb}{0.0, 0.62, 0.38}
\definecolor{darkgreen}{rgb}{0, 0.392, 0}
\definecolor{darkcyan}{RGB}{0,139,139}

\newtheorem{theorem}{Theorem}
\newtheorem{remark}[theorem]{Remark}

\newcommand{\bD}{\pmb{D}}

\newcommand{\br}{\pmb{r}}
\newcommand{\bv}{\pmb{v}}

\newcommand{\bu}{\pmb{u}}

\newcommand{\bn}{\pmb{n}}

\newcommand{\iS}{\mathcal{S}}
\newcommand{\iT}{\mathcal{T}}
\newcommand{\iG}{\mathcal{G}}
\newcommand{\iF}{\mathcal{F}}
\newcommand{\iN}{\mathcal{N}}

\newcommand{\iH}{\mathcal{H}}
\newcommand{\iB}{\mathcal{B}}
\newcommand{\Div}{\text{div}\;}

\newcommand{\prm}{\alpha}

\newcommand{\iR}{\mathcal{R}}
\newcommand{\iL}{\mathcal{L}}
\newcommand{\iLR}{\mathcal{L}^R_{\Omega}}
\newcommand{\iLD}{\iL_{\Omega}}
\newcommand{\iLDi}{\iL_{\Omega_i}}
\newcommand{\iD}{\mathcal{D}}
\newcommand{\iLRi}{\mathcal{L}^R_{\Omega_i}}
\newcommand{\mR}{\mathbb{R}}

\newcommand{\iK}{\mathcal{K}}

\newcommand{\iU}{\mathcal{U}}
\newcommand{\iE}{\mathcal{E}}

\newcommand{\ds}{\displaystyle}

\newcommand{\uOE}{u_{\Omega,E}}
\newcommand{\uiE}{u_{i,E}}
\newcommand{\uKE}{u_{K,E}}
\newcommand{\loh}{\Lambda_{h}}
\newcommand{\noh}{N_{h,\Omega}}

\newcommand{\lh}{\Lambda_{h}}
\newcommand{\iEhb}{\mathcal{G}_h}
\newcommand{\iEhzero}{\mathcal{G}_h^0}
\newcommand{\Dt}{\tau}
\newcommand{\lha}{\lambda_{h,a}}
\newcommand{\lhai}{\lambda_{h,a,i}}
\newcommand{\llh}{\lambda_{h}}
\newcommand{\llhi}{\lambda_{h,i}}

\usepackage[ruled]{algorithm}            
\usepackage{algorithmicx} 
\usepackage{algpseudocode}
\usepackage{subcaption}

\makeatletter
\renewcommand{\ALG@name}{Problem}
\makeatother

\algrenewcommand\algorithmicindent{1.0em}
\algrenewcommand{\algorithmicrequire}{\tiny \textbf{\tiny{Input:}}}
\algrenewcommand{\algorithmicensure}{\small{ \textbf{Output:}}}
\algrenewcommand{\alglinenumber}[1]{\color{red!80!blue}\footnotesize#1:}
\algrenewcommand{\alglinenumber}[1]{\footnotesize#1:}
\usepackage{calc} 
\algrenewcommand\algorithmicfunction{\widthof{\textbf{\small{Function}}}}
\algrenewcommand{\algorithmiccomment}[1]{{\color{Green}\hfill\footnotesize $\vartriangleright$ #1}} 
\algrenewcommand{\algorithmiccomment}[1]{{\hfill\footnotesize $\vartriangleright$ \textit{#1}}} 
\algnewcommand\algorithmicto{\textbf{to}}
 \algrenewtext{For}[3]%
   {\algorithmicfor\ $#1= #2$,..., $#3$ }
\makeatletter
\renewcommand{\ALG@beginalgorithmic}{\small}

\let\OldStatex\Statex
\renewcommand{\Statex}[1][3]{%
  \setlength\@tempdima{\algorithmicindent}%
  \OldStatex\hskip\dimexpr#1\@tempdima\relax}
\newcommand{\StatexIndent}[1][3]{%
  \setlength\@tempdima{\algorithmicindent}%
  \Statex\hskip\dimexpr#1\@tempdima\relax}
\makeatother

\newif\ifstartedinmathmode

\begin{document}

\begin{frontmatter}

\title{Space-time domain decomposition for advection-diffusion problems in mixed formulations}
\tnotetext[mytitlenote]{Partially supported by ANDRA, the French agency for nuclear waste management, and by GNR Momas.
}


\author[mymainaddress,mysecondaryaddress]{Thi-Thao-Phuong Hoang}
\author[mymainaddressLAGA,mymainaddress]{Caroline Japhet\corref{mycorrespondingauthor}}
\cortext[mycorrespondingauthor]{Corresponding author.}
\author[mymainaddress]{Michel Kern}
\author[mymainaddress]{Jean E. Roberts}


\address[mymainaddress]{INRIA Paris-Rocquencourt, 78153 Le Chesnay Cedex, France, Michel.Kern@inria.fr, Jean.Roberts@inria.fr}
\address[mymainaddressLAGA]{Universit\'e Paris 13, UMR 7539, LAGA,
  99 Avenue J-B Cl\'ement, 93430 Villetaneuse, France, japhet@math.univ-paris13.fr}
\address[mysecondaryaddress]{Ho Chi Minh City University of Pedagogy, Vietnam, hoang.t.thao.phuong@gmail.com.}

\begin{abstract}
This paper is concerned with the numerical solution of porous-media flow and transport problems, i. e. heterogeneous, advection-diffusion problems.  Its aim is to investigate numerical schemes for these problems in which different time steps can be used in different parts of the domain.  Global-in-time, non-overlapping domain-decomposition methods are coupled with operator splitting making possible the different treatment of the advection and diffusion terms.  Two domain-decomposition methods are considered: one uses the time-dependent Steklov–Poincar\'e operator and the other uses optimized Schwarz waveform relaxation (OSWR) based on Robin transmission conditions. For each method, a mixed formulation of an interface problem on the space-time interface is derived, and different time grids are employed to adapt to different time scales in the subdomains. A generalized Neumann-Neumann preconditioner is proposed for the first method. To illustrate the two methods numerical results for two-dimensional problems with strong heterogeneities are presented.  These include both academic problems and more realistic prototypes for simulations for the underground storage of nuclear waste.
\end{abstract}

\begin{keyword} 
   mixed formulations, domain decomposition, advection-diffusion,
  optimized Schwarz waveform relaxation,
time-dependent Steklov–Poincar\'e operator, nonconforming time grids
\end{keyword}

\end{frontmatter}


\section{Introduction}

In the simulation of contaminant transport in and around a nuclear waste repository the time scales vary over several orders of magnitude due to different material properties in different parts of the repository and to variations in the hydrogeological properties of the surrounding geological layers.  The domain of calculation is in fact a union of several regions with drastically different physical properties. In addition, different processes are taking place and these occur on very different time scales.  Thus the use of a single-size time-step throughout the domain and for all the physical processes involved is hardly reasonable.  For this reason we are interested in numerical schemes that allow different time steps for different parts of the domain and different time steps for different processes. For the different time steps for different parts of the domain we use global-in-time, non-overlapping domain-decomposition and for the different time steps for different processes we use operator splitting.  
For the space discretization, we use mixed finite elements~\cite{brezzi1991mixed,RobertsThomas} as they are mass conservative and they handle well heterogeneous and anisotropic diffusion tensors.

In \cite{PhuongSINUM}, two space-time domain-decomposition methods for the diffusion problem in a mixed setting were introduced. Both methods rely on a reformulation of the initial problem as a space-time interface problem, through the use of trace operators. The resulting interface problem can then be solved by various iterative methods.

For the first method, a global-in-time preconditioned Schur method (GTP-Schur), the trace operator in question is a time-dependent Dirichlet-to-Neumann (a.k.a. Steklov-Poincar\'e-type) operator; the interface problem is solved by preconditioned GMRES. 
Steklov-Poincar\'e-type methods are known to be very efficient for stationary problems with highly heterogeneous coefficients. See \cite{Achdou:DDnonsym:1999,Achdou:DDP:2000,Bourgat:1989,CowsarBDD,Mandelweights,Mathew:DDM:2008,quarteroni2008numerical,Toselli:DDM:2005} for more information concerning these methods for stationary problems.
Other extensions to time dependent diffusion problems are given in~\cite{gander2014dirichlet,Kwok}.
These are related to a relaxation method applied to the interface problem.

For the second method,  a global-in-time optimized Schwarz method (GTO-Schwarz), the associated trace operator is a time-dependent Robin-to-Robin-type (or possibly Ventcell-to-Ventcell-type) operator and the method is related to the optimized Schwarz waveform relaxation algorithm (OSWR).
These more general trace operators introduce additional coefficients that can be optimized to improve convergence rates, see~\cite{Bennequin,JeanRobinmixed,Gander2006,OSWRwave,HJKRDD22,JaphetDD9,Japhet-BIC-2001,VMartin}. Generalizations to heterogeneous problems and/or nonmatching time grids were introduced in~\cite{PMThesis,BertheDD21,OSWR2d,BlayoHJ,OSWR1d2,Haeberlein,Haeberlein-NSO-2013,OSWR3sub,OSWRDG,OSWRDG2,PhuongThesis,PhuongSINUM}. A suitable time projection to handle nonconforming grids in time between subdomains is defined using an optimal projection algorithm as in~\cite{Projection2d:3d,Projection1d} with no additional grid.
The OSWR method corresponds to solving the interface problem by Jacobi iterations. We will also use a GMRES method.

An extension of the GTP-Schur and GTO-Schwarz methods to a problem modeling the flow of a compressible fluid in a porous medium with a ``fast path'' fracture is given in~\cite{HJKRFracture}. More precisely, a reduced model is considered in which the fracture is treated as an interface between two subdomains and the associated trace operator is
of Ventcell-to-Robin-type.

This paper can be seen as a sequel to~\cite{PhuongSINUM}: we provide a new approach that
makes it possible to extend both the GTP-Schur and GTO-Schwarz methods to model coupled advection-diffusion problems. More precisely, because one may want to treat advection and diffusion with different numerical schemes, we use operator splitting within the subdomains~\cite{Arbogast:1996,hundsdorfer,Mazzia:Goudunov:2002,Mazzia:splitting:2000,Siegel:1997,Vasilevski:2008,Wheeler:1988}. This allows the use of different time steps from one subdomain to the next as well as within each subdomain, for the advection and the diffusion. For both the GTP-Schur and GTO-Schwarz methods, extensions of the discrete counterparts of the interface problems derived for the diffusion problem in~\cite{PhuongSINUM} are proposed.
New unknowns for the advection are introduced for the interface problems associated with Dirichlet transmission conditions between subdomains for the advection step.
For GTP-Schur, an extension of the Neumann-Neumann preconditioner of~\cite{PhuongSINUM} to the advection-diffusion problem is introduced.
Advection is approximated in time with the explicit Euler method and in space with an upwind, cell-centered finite volume method, while diffusion is approximated in time with the implicit Euler method and in space with a mixed finite element method. It has been shown that treating the advection explicitly can significantly reduce the numerical diffusion (see e.g.~\cite{ChertockKurganov}). 
An upwind operator is introduced to simplify the formulation of the advection step.
This operator is useful for extending the discrete counterparts of the interface problems derived for the diffusion problem in~\cite{PhuongSINUM}.

\bigskip

This paper is organized as follows: in the next section we present both the model problem in mixed
form and the fully discrete problem. In Section~\ref{A2Sec:DD}
we derive the discrete multidomain problem in an operator splitting context and define two discrete interface problems, extensions of the GTP-Schur and GTO-Schwarz methods analyzed in~\cite{PhuongThesis,PhuongSINUM}. For the GTP-Schur method, a generalized Neumann-Neumann preconditioner is given. 
In Section~\ref{Sec:AdvTime} we describe how we handle the nonconforming time grids (for advection and diffusion time steps) using $ L^{2} $ projections.
In Section~\ref{A2Sec:Num}, results of 2D numerical experiments, both academic experiments and more realistic prototypes for nuclear waste disposal simulation, showing that the methods preserve the order of the global scheme are presented.
The behaviors of the two methods are discussed and compared.

\section{Operator splitting for a model problem in a single domain}
\label{A2Sec:Model}
For a bounded domain $ \Omega $ of $ \mathbb{R}^{d} \; (d=2,3) $ with Lipschitz boundary $ \partial \Omega $ and some fixed time $T>0$, consider the linear advection-diffusion problem written in mixed form:
\begin{equation} \label{A2model}
\begin{array}{rll} \phi \partial_{t} c + \Div(\bu c+ \br) & =f & \text{in} \; \Omega \times (0,T ),\\
\nabla c + \bD^{-1} \br &=0  & \text{in} \; \Omega \times (0,T), \\
c&=0 & \text{on} \; \partial \Omega \backslash \Gamma \times (0,T), \\
c&=c_D & \text{on} \; \Gamma \times (0,T), \\
c(\cdot, 0) & = c_{0} & \text{in} \; \Omega,
\end{array} 
\end{equation}
where $ c $ is the concentration of a contaminant dissolved in a fluid, $ f $ the source term, 
$ \phi $ the porosity, $ \bu $ the Darcy velocity (assumed to be given and time-independent), $\bD$ a symmetric time independent diffusion tensor, $\Gamma$ a part of the boundary and $c_D$ Dirichlet boundary data on $\Gamma$.
We have singled out a part of the boundary for non-homogeneous Dirichlet data because it will be useful when we derive the GTP-Schur method in
  Section~\ref{A2Sec:DD}. The monodomain or global problem corresponds to the case $\Gamma=\emptyset$.
  For simplicity we have imposed only Dirichlet boundary conditions. The analysis presented in the following can be generalized to other types of boundary conditions.

For the time discretization, we use a splitting method, for solving problem~\eqref{A2model}: the advection equation is approximated by the forward Euler method and the diffusion equation by the backward Euler method. The resulting scheme is first-order accurate in time, $ O \left (\Dt\right ) $ (see, e.g., \cite{Dawson93,Leveque:1992}).
We consider a locally mass-conserving approximation scheme, more specifically an upwind, cell-centered finite volume method for the advection equation, and a mixed finite element method for the diffusion equation. 
Below, we give the fully discrete problem associated with these discretization techniques. 

Let $ \iK_{h} $ be a finite element partition of $ \Omega $ into rectangles and let $\iE_h$ be the set of all faces
of elements of $ \iK_{h} $, $\iEhb$ the set of those lying on $\Gamma$ and $\iEhzero$ those in the interior.
For $K \in \iK_{h} $, let $\bn_K$ be the unit, normal, outward-pointing vector field on $\partial K$ and let $\bn_\Omega$ be the unit, normal, outward-pointing vector field on $\partial \Omega$.
For simplicity, we suppose $ \Omega \subset \mR^{2}$. We use the lowest order Raviart-Thomas 
mixed finite element spaces $ M_{h} \times \Sigma_{h}  \subset L^{2}(\Omega) \times H(\text{div}, \Omega)$ (see, e.g., \cite{brezzi1991mixed, RobertsThomas}) : 
\begin{align*}
M_{h} & = \left \{ \mu \in L^{2}(\Omega): \mu_{\mid K} = \text{const}, \, K \in \iK_{h} \right \}, \\
\Sigma_{h} & = \left \{ \bv \in H(\text{div}, \Omega): \bv_{\mid K} = \left (a_{K} + b_{K} x, c_{K} +d_{K} y \right ), \, (a_{K}, b_{K}, c_{K}, d_{K}) \in \mR^{4}, \; K \in \iK_{h} \right \}. 
\end{align*}
The degrees of freedom of $ c_{h} \in M_{h} $ correspond to the average values of $ c_{h} $ on the elements $ K \in \iK_{h} $, and those of $ \br_{h} $ correspond to the values of the flux of $ \br_{h}$ across the edges $ E $ of $ K $.

We shall also make use of the spaces
\begin{align}
  \loh &=\{\lambda \in L^2(\Gamma) \, : \, \lambda_{|E}=\text{const}, \ E \in \iEhb\} \label{eq:loh-1}\\
  \noh &=\{\eta \in L^2(\bigcup_{E\in \iE_h} E) \, : \, \eta_{|E}=\text{const}, \ E \in \iE_h\}\label{eq:loh-2}
  \end{align}
for the approximation of the boundary data and for the upwind values respectively.
For the time grid, we consider, for simplicity, a uniform partition of $ \left (0,T\right ) $ into $ N $ subintervals $ \left (t^{n},t^{n+1}\right ) $ with length $ \Dt=t^{n+1}-t^{n} $ for $ n~=~0, \hdots, N-1 $ with $ t^{0} = 0 $ and $ t^{N} = T $ (the derivation can be easily generalized to the case of nonuniform partitions). 
In order to satisfy the CFL condition required for the explicit scheme used for the advection equation without imposing that condition on the diffusion equation, we consider sub-time steps for the advection part: $ \Dt_{a}  = \Dt /L, \; \text{for some } \; L \geq 1, $ and $ t^{n,l} = l\Dt_{a} + t^{n}, \; \text{for} \; \; l=0, \hdots, L$,
and $n=0, \hdots, N-1. $
Note that $ t^{n,0} = t^{n} $ and $ t^{n,L} = t^{n+1} $. 

The operator splitting algorithm is initialized by defining $ c_{h}^{0} $ to be the $L^{2}$ projection of $ c_{0} $
onto~$ M_{h} $:
\begin{equation}\label{eq:ch0}
c_{h \mid {K}}^{0}:= \frac{1}{|K|} \ds\int_{K} c_{0}, \; \forall K \in \iK_{h},
\end{equation}
and $\lha^{n,l}, \ 0\le l \le L-1, \ 0\le n \le N,$ to be the $L^{2}$ projection of $c_D(\cdot,t^{n,l})$ onto $ \loh $:
\begin{equation}\label{eq:lambda_n_l}
(\lha^{n,l})_{|E}:= \frac{1}{|E|} \ds\int_{E} c_D(\cdot,t^{n,l}), \; \forall E \in \iEhb, 
\end{equation}
where $|K|$, respectively $|E|$, denotes the measure of $K$, respectively $E$.
For convenience of notation, we also write $ c_{h}^{0,0} $ for $ c_{h}^{0} $ and $ \lha^{0,0} $, for
$ \lha^{0} $. We will find it useful to use different notation for the boundary values that we will use for the diffusion step and the boundary values that we will use for the advection step: $\lha^{n,l}$ is defined
in~\eqref{eq:lambda_n_l} and is used for the advection step whereas
$\llh^n:=\lha^{n,0}$
is used for the diffusion step.
\\
For $ n=0, \hdots, N-1$, at step $n$, we first compute $ c_{h}^{n,l} $, the approximation of $ c(t^{n,l}) $, for $l=1,  \hdots, L$ using the advection equation and then we compute $ c_{h}^{n+1} $ and $ \br_{h}^{n+1} $, approximations of $ c(t^{n+1}) $ and $ \br (t^{n+1}) $ respectively, using the diffusion equation.  
\\
As we use an upwind scheme for the advection equation, to compute $ c_{h}^{n,l+1} $ for $ n=0, \hdots,~N-1$, $l=0, \hdots, L-1, $ in addition to the value $ c_{h}^{n,l} $, we will need an upwind value $ \hat{c}_{h}^{n,l} $ of the concentration on each edge of the grid.
This value depends on the Darcy velocity $\bu$ and we define for each edge~$E \in \iE_h$ 
\begin{equation*} \label{eq:average}
  \uOE=\frac{1}{|E|}\int_E \bu \cdot \bn_{\Omega}, \qquad
  \uKE=\frac{1}{|E|}\int_E \bu \cdot \bn_{K}, \ \forall K\in \iK_{h}.
\end{equation*}
The upwind value is defined with an upwind operator $ \iU_{h}$ that associates to an element in
$M_{h} \times \loh$, a value in $ N_{h} $. These latter values will be the upwind values:
\begin{eqnarray} 
\iU_{h}: \; M_{h} \times \loh \rightarrow N_{h} \hspace{5cm}\nonumber\\[2mm]
\left (\iU_{h}(c_{h},\lha)\right )_{\mid_{E}} =
\begin{cases}
  (\lha)_{|E} & \text{if } E \in \iEhb \text{ and } \uOE<0
  \text{ (fluid entering } \Omega  \text{ through } E \text{)}, \\
  0 & \text{if } E \in \iEhb \text{ and } \uOE \ge 0, \\
  (c_{h})_{|K} & \text{if } E \in \iEhzero \text{ is an edge of } K
  \text{ and  } \uKE\ge 0 \text{ (fluid exiting } K \text{ through } E
    \text{)},  \\
    0 & \text{if } E\in \iEhzero \text{ and } \uKE = 0 \text{ for some
    } K \in \iK_h. \label{A2UpwindO}
\end{cases}
\end{eqnarray}
The discrete problem with operator splitting is then defined as follows: 
\begin{algorithm}[H]
\captionsetup{font=footnotesize}
\caption{Discrete problem with operator splitting and Dirichlet boundary conditions}
\label{algo:splittingD}
\begin{algorithmic}[0]
\For{n}{0}{N-1}
   \State define $ c_{h}^{n,0} =  c_{h}^{n} $, where $ c_h^{0}$ is defined in \eqref{eq:ch0},
 \For{l}{0}{L-1}
    \State 1. define the upwind value $\hat{c}_{h}^{n,l}:=\iU_{h}(c_{h}^{n,l},\lha^{n,l})$
                  using \eqref{eq:lambda_n_l}, 
     \State 2. solve the advection equation
			\begin{equation} \label{A2time-advection}
			\ds\int_{K} \phi \frac{c_{h}^{n,l+1} - c_{h}^{n,l} }{\Dt_{a}}
                        + \sum_{E \subset \partial K} \ds \left(\hat{c}^{n,l}_{h}\right )_{\mid E} |E| \uKE =0,
                        \quad  \forall K \in \iK_{h},
			\end{equation} 
                 \State  to~obtain~$ c_{h}^{n,l+1} $.
                  The solution generated after these $ L $ advection steps is~$ c_{h}^{n,L} $. 
      \EndFor
   \State solve the diffusion equation
   \begin{equation} \label{A2time-diffusion}
			\begin{array}{rll} \ds\int_{K} \phi \frac{c_{h}^{n+1}-c_{h}^{n, L} }{\Dt} + \ds\int_{K} \Div \br_{h}^{n+1} & = \ds\int_{K} f(t^{n+1}), & \forall K \in \iK_{h} , \vspace{0.15cm}\\
			\ds\int_{\Omega} \bD^{-1} \br_{h}^{n+1} \cdot \bv
                        -\ds\int_{\Omega}  c_{h}^{n+1} \Div \bv
                        & = \int_\Gamma \llh^{n+1} (\bv \cdot \bn), & \forall \bv \in \Sigma_{h},  
			\end{array}
			\end{equation}
  \State to obtain $ c_{h}^{n+1} $ and $ \br_{h}^{n+1} $.
  \EndFor
\end{algorithmic}
\end{algorithm}

(Recall that the test functions for the first equation in~\eqref{A2time-diffusion} are just linear combinations of the characteristic functions on the elements $ K \in \iK_{h} $.)
%
%
%
  Using the notation
\begin{align*}
  \lambda_{a} = \left (\lha^{n,l}\right )_{n=0, \hdots, N-1, \; l=0, \hdots, L-1}, & \qquad \qquad \lambda = \left ( \llh^{n+1} \right )_{n=0, \hdots, N-1},\\
 \hat{c}_{h}^{\Dt, \Dt_{a}} = \left (\hat{c}_{h}^{n,l}\right )_{n=0, \hdots, N-1, \; l=0, \hdots, L-1},  & \qquad \qquad \left (c_{h}^{\Dt}, \br_{h}^{\Dt}\right ) = \left (c_{h}^{n}, \br_{h}^{n}\right )_{n=1, \hdots, N},
 \end{align*}
we can define a discrete solution operator $ \iLD$ as follows
\begin{equation*} 
\hspace{-0.3cm}\begin{array}{rl} \iLD: \loh^{N \times L} \times \loh^{N} \times L^{2}(0,T; L^{2}(\Omega)) \times H_{*}^{1}(\Omega) & \rightarrow  \left (N_{h}\right) ^{N\times L}  \times \left (M_{h}\right) ^{N} \times \left (\Sigma_{h}\right) ^{N} \vspace{0.15cm}\\
\left (\lambda_{a}, \lambda, f, c_{0}\right )  &\mapsto \left (\hat{c}_{h}^{\Dt, \Dt_{a}}, c_{h}^{\Dt}, \br_{h}^{\Dt}\right ), 
\end{array} \vspace{0.15cm}
\end{equation*} 
where $\left (\hat{c}_{h}^{\Dt, \Dt_{a}}, c_{h}^{\Dt}, \br_{h}^{\Dt}\right )$ is the solution
of Problem~\ref{algo:splittingD}, and $H_{*}^{1}(\Omega)=\{v \in H^1(\Omega) ; v=0 \text{ on } \partial \Omega \backslash  \Gamma\}$.
We remark that as the advection is approximated explicitly, in the definition of $ \iLD$ we have extracted the upwind values $ \hat{c}_{h}^{n,l} $ for $ l=0, \hdots, L-1 $ (instead of $ l=1, \hdots, L $) for each $ n, \; n=0, \hdots, N-1 $. 

In the next section, we consider the domain decomposition approach for solving Problem~\ref{algo:splittingD}. An equivalent multidomain problem adapted to the splitting approach will be formulated and from that we will derive two global-in-time domain decomposition methods.
%
%
%
\section{Domain decomposition with operator splitting}
\label{A2Sec:DD}
For simplicity, we consider a decomposition of a domain $ \Omega $ into just two non-overlapping subdomains $ \Omega_{1} $ and $ \Omega_{2} $ (the analysis can be generalized to the case of any number of subdomains). We assume that the partitions $ \iK_{h,1} $ of subdomain $ \Omega_{1} $ and $ \iK_{h,2} $ of subdomain $ \Omega_{2} $ are such that their
union $ \iK_{h}= \bigcup_{i=1}^2 \iK_{h,i}  $ forms a finite element partition of $ \Omega $ . We denote by $ \Gamma:=\partial \Omega_{1} \cap  \partial \Omega_{2} $ the interface between the subdomains, and denote by $ \iG_{h} $ the set of edges (or faces) of elements of $ \iK_{h} $ that lie on $ \Gamma $.
For $i=1,2$, let $ \bn_{i} $ denote the unit, normal, outward-pointing vector field on $\partial\Omega_i$, and for any scalar, vector or tensor valued function $ \psi $ defined on $\Omega$, let $\psi_i$ denote the restriction  of $ \psi$ to $ \Omega_{i} $. 

\noindent
We define the set of the inflow and outflow boundary edges on the interface for each subdomain:
$$ \iG_{h,i}^{\text{in}}~:=~\left \{ E \in \iG_{h}: \ds\uiE <0 \right \},
\quad \iG_{h,i}^{\text{out}}~:=\iG_{h,j}^{\text{in}}
\qquad \text{for} \; i=1,2, \; j=3-i,$$
where
$
  \uiE=\frac{1}{|E|}\int_E \bu \cdot \bn_{i}.
$
Thus
$\iG_{h} \setminus (\iG_{h,1}^{\text{in}} \cup \iG_{h,2}^{\text{in}}) = \{E \in \iE_h; \; \uiE = 0\}$,
and $\iG_{h,1}^{\text{in}} \cap \iG_{h,2}^{\text{in}} = \emptyset$.
%
%

\begin{figure}[H]
\vspace{0.3cm}
\hspace{1cm} \begin{minipage}[c]{0.5 \linewidth}
\setlength{\unitlength}{1pt} 
\begin{picture}(140,170)(0,0)
\thicklines
\put(0,0){\includegraphics[scale=0.81]{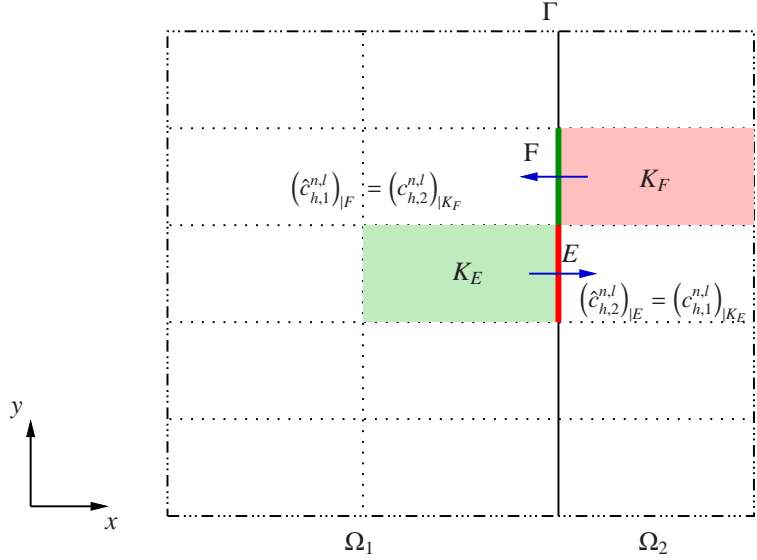} \\}
\put(194,188){$ \Gamma $}
\put(30,-4){$x$}
\put(-5,40){$y$}
\put(120,-12){$\Omega_{1}$}
\put(230,-12){$ \Omega_{2} $}
\put(201,97){$ E $}
\put(160,90){$K_{E}$}
\put(208,80){\small $ \left (\hat{c}_{h,2}^{n,l}\right )_{\mid E }  = \left (c_{h,1}^{n,l}\right )_{\mid K_{E}}$}
\put(230,125){$ K_{F} $}
\put(187,135){F}
\put(100,122){\small $ \left (\hat{c}_{h,1}^{n,l}\right )_{\mid F}  \,= \left (c_{h,2}^{n,l}\right )_{\mid K_{F}}$}
\end{picture}
\end{minipage} \vspace{0.5cm}
\caption{An illustration of the upwind concentration in the context of domain decomposition, the arrows represent the direction of the normal flux across the edges, $ \ds\int_{E} \bu \cdot \bn, \; $ (for a fixed normal vector $ \bn = (1,0) $).}
\label{A2Fig:UpwindMultidom} 
\end{figure}

Let $ M_{h} $ and $ \Sigma_{h} $ denote the mixed finite element spaces as defined in Section~\ref{A2Sec:Model}, and let
$ M_{h,i} $ and $ \Sigma_{h,i} $, $ i=1,2, $ be the spaces of restrictions of the functions in these spaces
to $ \Omega_{i} $. To define the transmission conditions, we will use the space $\loh$ defined in~\eqref{eq:loh-1}
(where $\Gamma$ is the interface).
Because the discretization in space is conforming, we have
$$ \lh=\{(\bv\cdot \bn_i)_{\mid_{\Gamma}} : \; \bv \in \Sigma_{h,i}\}, \ i=1,2. 
$$
Now, it is  straightforward to show that the global discrete problem
  defined in Problem~\ref{algo:splittingD} with $c_D=0$ (i.e. $\lha=\llh=0$)
is equivalent to the following multidomain problem:
\begin{algorithm}[H]
\captionsetup{font=footnotesize}
\caption{Multidomain problem with operator splitting and physical transmission conditions}
\label{algo:splittingmulti}
\begin{algorithmic}[0]
\State Solve in each subdomain $\Omega_i$, i=1,2:
  \For{n}{0}{N-1}
   \State define $ c_{h,i}^{n,0} = c_{h,i}^{n} $, where $ \left (c_{h,i}^{0}\right )_{\mid K}:= \frac{1}{|K|} \ds\int_{K} c_{0}, \; \forall K \in \iK_{h,i}, $
 \For{l}{0}{L-1}
 	\State 1. define $\left (\lhai^{n,l}\right )_{\mid E}$ for $E \in \iG_{h,i}^{\text{in}}$ through
         the transmission condition
        \begin{equation} \label{A2TCs-adv}
\left (\lhai^{n,l}\right )_{\mid E} = \left (c_{h,j}^{n,l}\right )_{\mid K_{E}}, \; 
        \forall E \in \iG_{h,i}^{\text{in}}, \quad \text{with } \, j=3-i,
\end{equation}
        \State where $ K_{E} $ is the element that has $ E $ as an edge. Then
        define the upwind values
        $ \hat{c}_{h,i}^{n,l} = \iU_{h} \left (c_{h,i}^{n,l}, \lhai^{n,l}\right )$.
        \State 2. solve the advection equation~\eqref{A2time-advection} in $\Omega_i$:
	\begin{equation} \label{A2time-advection-M1_i}
	  \ds\int_{K} \phi \frac{c_{h,i}^{n,l+1} - c_{h,i}^{n,l} }{\Dt_{a,i}}
          + \sum_{E \subset \partial K} \ds \left(\hat{c}^{n,l}_{h,i}\right )_{\mid E} |E| \uKE =0,
\quad \forall K \in \iK_{h,i}, 
	\end{equation}
        \State to obtain $ c_{h,i}^{n,l+1} $. The solution generated after these $ L $ advection steps is~$ c_{h,i}^{n,L} $. 
      \EndFor
   \State solve the diffusion equation~\eqref{A2time-diffusion} in $\Omega_i$:
     \begin{equation} \label{A2time-diffusion-M1-i}
     \hspace{-0.7cm}
     \begin{array}{rll} \ds\int_{K} \phi \frac{c_{h,i}^{n+1}-c_{h,i}^{n, L} }{\Dt_i}
       + \ds\int_{K} \Div \br_{h,i}^{n+1} & = \ds\int_{K} f(t^{n+1}), & \forall K \in \iK_{h,i} , \\
       \ds\int_{\Omega_i} \bD^{-1} \br_{h,i}^{n+1} \cdot \bv -\ds\int_{\Omega_i}  c_{h,i}^{n+1} \Div \bv
          & =  \ds\int_{\Gamma} \llhi^{n+1} (\bv \cdot \bn_i), & \forall \bv \in \Sigma_{h,i},
     \end{array}
   \end{equation}
    \State together with the transmission conditions
\begin{equation} \label{A2TCs-diff}
  \ds\int_{E} \lambda^{n+1}_{h,1} = \ds\int_{E} \lambda^{n+1}_{h,2}, \qquad 
\ds\int_{E} \left (\br_{h,1}^{n+1} \cdot \bn_{1} +\br_{h,2}^{n+1} \cdot \bn_{2} \right )= 0,
\quad \forall E \in \iG_{h},  
\end{equation}
 \State to obtain $ c_{h,i}^{n+1} $ and $ \br_{h,i}^{n+1} $.  
  \EndFor
\end{algorithmic}
\end{algorithm}
Equation~\eqref{A2TCs-adv} serves as a Dirichlet boundary condition on $ E \in \iG_{h,i}^{\text{in}} $ and it defines the transmission condition for the advection equation. As for the pure diffusion problem~\cite{PhuongSINUM}, the transmission conditions~\eqref{A2TCs-diff} for the diffusion equation consist of the equality between the concentration and the conservation of the normal diffusive flux across the interface.
\begin{remark}
Because of the decomposition, the upwind value of $ c_{h,i}^{n,l} $ on edges on the interface $ \Gamma $ may not depend only on the element values $ c_{h,i}^{n,l} $ inside the subdomain $ \Omega_{i} $. In particular, if $ \iG_{h,i}^{\text{in}} \neq \emptyset $ (i.e. there is fluid flowing into $ \Omega_{i} $ through some part of $\, \Gamma $), the transmission condition~\eqref{A2TCs-adv} means that the upwind concentration $ \hat{c}_{h,i}^{n,l} $ on the edge $ E \in \iG_{h,i}^{\text{in}}$ is defined by the concentration value of the element in the neighboring subdomain (see Figure~\ref{A2Fig:UpwindMultidom}):		
\begin{equation} \label{A2TCs-adv-2}
\left (\hat{c}_{h,i}^{n,l}\right )_{\mid E} = \left (c_{h,j}^{n,l}\right )_{\mid K_{E}}, \; \forall E \in \iG_{h,i}^{\text{in}}, 
\end{equation}
where $ K_{E} $ necessarily in $\iK_{h,j}, \; j=3-i,$ is the element that has $ E $ as an edge.
\end{remark}
Alternatively, and equivalently to \eqref{A2TCs-diff}, one may impose Robin transmission conditions (for the diffusion equation), for all  $E \in \iG_{h} $ and $ n=0, \hdots, N-1:$
\begin{equation} \label{A2RobinTCs-diff}
\begin{array}{ll} \ds\int_{E}  \left (-\br_{h,1}^{n+1} \cdot \bn_{1} + \prm_{1,2} \lambda_{h,1}^{n+1}\right ) & = \ds\int_{E}  \left (\br_{h,2}^{n+1} \cdot \bn_{2} + \prm_{1,2} \lambda_{h,2}^{n+1}\right ),\vspace{0.15cm}\\
\ds\int_{E} \left (-\br_{h,2}^{n+1} \cdot \bn_{2} + \prm_{2,1} \lambda_{h,2}^{n+1} \right )& = \ds\int_{E} \left (\br_{h,1}^{n+1} \cdot \bn_{1}+ \prm_{2,1} \lambda_{h,1}^{n+1}\right ), \vspace{0.15cm}
\end{array} 
\end{equation}
where $\prm_{1,2}$ and $\prm_{2,1}$ are two positive constants.
The GTP-Schur method is based on
\eqref{A2time-advection-M1_i}, \eqref{A2time-diffusion-M1-i} together with the "natural" transmission conditions~\eqref{A2TCs-adv}, \eqref{A2TCs-diff} while the GTO-Schwarz method is based on~\eqref{A2time-advection-M1_i}, \eqref{A2time-diffusion-M1-i} together with transmission conditions \eqref{A2TCs-adv}, \eqref{A2RobinTCs-diff}. For both methods the multidomain problem is formulated as a problem posed on the space-time interface via the use of interface operators.
These interface operators are defined via the solution of local Robin problems in the subdomains.

For a given Robin parameter $\prm$, we can define a discrete Robin solution operator $ \iLR$ as follows
\begin{equation}\label{eq:LR}
\hspace{-0.3cm}\begin{array}{rl} \iLR: \loh^{N \times L} \times \loh^{N} \times L^{2}(0,T; L^{2}(\Omega)) \times H_{*}^{1}(\Omega) & \rightarrow  \left (N_{h}\right) ^{N\times L}  \times \left (M_{h}\right) ^{N} \times \left (\Sigma_{h}\right) ^{N} \vspace{0.15cm}\\
\left (\lambda_{a}, \xi, f, c_{0}\right )  &\mapsto \left (\hat{c}_{h}^{\Dt, \Dt_{a}}, c_{h}^{\Dt}, \br_{h}^{\Dt}\right ), 
\end{array} \vspace{0.15cm}
\end{equation} 
where $\left (\hat{c}_{h}^{\Dt, \Dt_{a}}, c_{h}^{\Dt}, \br_{h}^{\Dt}\right )$ is the solution of Problem~\ref{algo:splittingR}:
\begin{algorithm}[H]
\captionsetup{font=footnotesize}
\caption{Discrete local problem with operator splitting and Robin boundary conditions}
\label{algo:splittingR}
\begin{algorithmic}[0]
\For{n}{0}{N-1}
   \State define $ c_{h}^{n,0} =  c_{h}^{n} $, where $ c_h^{0}$ is defined in \eqref{eq:ch0},
 \For{l}{0}{L-1}
 	\State 1. define the upwind value $\hat{c}_{h}^{n,l}:=\iU_{h}(c_{h}^{n,l},\lha^{n,l})$
                  using \eqref{eq:lambda_n_l}, 
        \State 2.  solve the advection equation
			\begin{equation} \label{A2time-advection-R}
			\ds\int_{K} \phi \frac{c_{h}^{n,l+1} - c_{h}^{n,l} }{\Dt_{a}}
                        + \sum_{E \subset \partial K} \ds \left(\hat{c}^{n,l}_{h}\right )_{\mid E} |E| \uKE =0,
                        \quad  \forall K \in \iK_{h},
			\end{equation} 
              \State    to~obtain~$ c_{h}^{n,l+1} $.
                  The solution generated after these $ L $ advection steps is~$ c_{h}^{n,L} $. 
      \EndFor
   \State solve the diffusion equation
\begin{equation} \label{A2time-diffusion-R}
	   \begin{array}{rll} \ds\int_{K} \phi \frac{c_{h}^{n+1}-c_{h}^{n, L} }{\Dt}
             + \ds\int_{K} \Div \br_{h}^{n+1}  = \ds\int_{K} f(t^{n+1}), & \forall K \in \iK_{h} ,&
             \vspace{0.15cm}\\
	     \ds\int_{\Omega} \bD^{-1} \br_{h}^{n+1} \cdot \bv -\ds\int_{\Omega}  c_{h}^{n+1} \Div \bv
             + \frac{1}{\prm}\ds\int_{\Gamma} (\br_{h}^{n+1} \cdot \bn) (\bv \cdot \bn)&= & \\
	      & \hspace{-2cm} - \frac{1}{\prm}\ds\int_{\Gamma} \xi_{h}^{n+1} (\bv \cdot \bn),
              & \hspace{-1cm}\forall \bv \in \Sigma_{h},
         \end{array}
         \end{equation}
 \State to obtain $ c_{h}^{n+1} $ and $ \br_{h}^{n+1} $.
  \EndFor
\end{algorithmic}
\end{algorithm}
\noindent
Note that the boundary terms in~\eqref{A2time-diffusion-R} correspond to the boundary condition
\begin{equation}\label{eq:RobinContinu}
  \alpha c_h^{n+1}-\br_h^{n+1} \cdot \bn=\xi_{h}^{n+1} \ \text{ on } \Gamma \times (0,T).
\end{equation}
\begin{remark}
As pointed out above, with the upwind scheme, the solution of the advection equation in a subdomain depends not only on the information in the subdomain and on its boundary, but also on information coming from the neighboring subdomain, while for the diffusion equation the solution is local to the subdomain as in the pure diffusion case. However, since we use operator splitting, we do not have the problem of slow convergence of the OSWR algorithm that has been observed when a fully implicit scheme and an upwind scheme for the advection are used \cite{Haeberlein}. With operator splitting we obtain separate transmission conditions for the advection part and for the diffusion part. In fact, we observe numerically (see Section~\ref{A2Sec:Num}) that the convergence is governed by the Robin transmission conditions associated with the diffusion equation and the optimized Robin parameters significantly improve the convergence of the algorithm (for both advection-dominated and diffusion-dominated problems). In our observations, the advection plays little or no role in the rate of convergence. 
\end{remark}
\begin{remark} \label{A2rmrkDiffStepAdvstep2M}
As the advection and the diffusion equations are treated separately, the formulations of the diffusion equation corresponding to the GTP-Schur method and the GTO-Schwarz method will be derived just as in~\cite{PhuongSINUM}. The formulation for the advection equation will be the same for both methods. 
\end{remark}
In the following, using operator splitting, we derive discrete interface problems for the advection-diffusion equation~\eqref{A2model} (with $c_D=0$), which are extensions of the discrete counterparts of the interface problems derived for the diffusion problem in~\cite{PhuongSINUM}.
%
%
%
%
\subsection{Global-in-time preconditioned Schur (GTP-Schur): an extension of the time-dependent Steklov-Poincar\'e operator approach}\label{sec:GTP-Schur}
Let $H_{*}^{1} (\Omega_i)=\{v\in H^{1} (\Omega_i), \; v=0 \text{ over } \partial \Omega_i \cap \partial \Omega \}, \, \text{for} \, i=1,2.$
%
We introduce the solution operators $ \iD_{i}:=\iLDi,\; i~=~1,2, $ which associates to an $ L^{2}(0,T; L^{2}(\Omega_{i})) $ source term $ f $ together with $ H_{*}^{1} (\Omega_{i}) $ initial data $ c_{0} $ and discrete boundary data $ \left (\lambda_{a}, \lambda\right ) $ given on $ \Gamma \times (0,T)$, the solution of Problem~\ref{algo:splittingD}
in $ \Omega_{i} \times (0,T) $ :
\begin{equation*} 
\hspace{-0.3cm}\begin{array}{rl} \iD_{i}: \lh^{N \times L} \times \lh^{N} \times L^{2}(0,T; L^{2}(\Omega_{i})) \times H_{*}^{1}(\Omega_{i}) & \rightarrow  \left (N_{h,i}\right) ^{N\times L}  \times \left (M_{h,i}\right) ^{N} \times \left (\Sigma_{h,i}\right) ^{N} \vspace{0.15cm}\\
\left (\lambda_{a}, \lambda, f, c_{0}\right )  &\mapsto
\left (\hat{c}_{h,i}^{\Dt, \Dt_{a}}, c_{h,i}^{\Dt}, \br_{h,i}^{\Dt}\right )=\iLDi\left(\lambda_{a}, \lambda, f, c_{0}\right ). 
\end{array} \vspace{0.15cm}
\end{equation*} 
For the problem on the interface, we will need as input from the subdomain problems the first component $\hat{c}_{h,i}^{\Dt, \Dt_{a}} $ of the output of $ \iD_{i} $ (for the advection step) and the values of the third component $ \br_{h,i}^{\Dt} $ (for the diffusion step). In fact, we need only the values of $\hat{c}_{h,i}^{\Dt, \Dt_{a}} $ associated with edges $ E $ in $ \iG_{h,i}^{\text{in}} $ and values of $ \br_{h,i}^{\Dt} $ associated with edges $ E $ in $ \iG_{h} $. Thus we define the two projection operators $ \iH_{i} $ and $ \iF_{i} $ as follows \vspace{0.15cm}
\begin{equation*}
\begin{array}{rl}  \iH_{i}: \left (N_{h,i}\right) ^{N\times L}  \times \left (M_{h,i}\right) ^{N} \times \left (\Sigma_{h,i}\right) ^{N} & \rightarrow  \left (\lh\right) ^{N\times L}    \vspace{0.15cm}\\
\left (\hat{c}_{h,i}^{\Dt, \Dt_{a}}, c_{h,i}^{\Dt}, \br_{h,i}^{\Dt}\right )  & \mapsto \left \{ \begin{array}{cl} 
0, & \forall E \in \iG_{h,i}^{\text{in}}, \vspace{0.15cm}\\
\left (\hat{c}_{h,i}^{\Dt, \Dt_{a}}\right )_{\mid {E}} , & \forall E \in \iG_{h,i}^{\text{out}},  
\end{array} \right .
\end{array} \vspace{0.1cm}
\end{equation*}
and \vspace{0.1cm}
\begin{equation*}
\begin{array}{rl} \iF_{i}: \left (N_{h,i}\right) ^{N\times L}  \times \left (M_{h,i}\right) ^{N} \times \left (\Sigma_{h,i}\right) ^{N} & \rightarrow \left (\lh\right) ^{N}    \vspace{0.15cm}\\
\left (\hat{c}_{h,i}^{\Dt, \Dt_{a}}, c_{h,i}^{\Dt}, \br_{h,i}^{\Dt}\right ) & \mapsto \left (\br_{h,i}^{\Dt} \cdot \bn_{i}\right )_{\mid {E}}, \;\forall E \in \iG_{h}.
\end{array} \vspace{0.2cm}
\end{equation*}
With these operators, we can rewrite the transmission condition~\eqref{A2TCs-adv} for the advection equation equivalently as 
\begin{equation*} 
\begin{array}{lll}
\ds\int_{t^{n,l}}^{t^{n,l+1}} \ds\int_{E}\lambda_{a} -  \iH_{1}\iD_{1}(\lambda_{a}, \lambda, f,c_{0}) & = 0, & \forall E \in \iG_{h,2}^{\text{in}},\vspace{0.15cm}\\
\ds\int_{t^{n,l}}^{t^{n,l+1}} \ds\int_{E}\lambda_{a} -  \iH_{2}\iD_{2}(\lambda_{a}, \lambda, f,c_{0}) & = 0, & \forall E \in \iG_{h,1}^{\text{in}}, \vspace{0.15cm}\\
& & \hspace{-4.5cm}\forall n=0, \hdots, N-1, \; \forall l=0, \hdots, L-1,
\end{array}
\end{equation*}
or 
\begin{equation} \label{A2IF-adv-M1}
\begin{array}{c}
\ds\int_{t^{n,l}}^{t^{n,l+1}} \ds\int_{E}\lambda_{a} -  \iH_{1}\iD_{1}(\lambda_{a}, \lambda, f,c_{0}) -  \iH_{2}\iD_{2}(\lambda_{a}, \lambda, f,c_{0}) = 0,\vspace{0.15cm}\\
\hspace{3cm} \forall E \in \iG_{h}, \; \; \forall n=0, \hdots, N-1, \; \forall l=0, \hdots, L-1. 
\end{array}
\end{equation}
Since we have imposed a Dirichlet condition on $ \Gamma $ 
for the diffusion equation, the first equation of~\eqref{A2TCs-diff} is satisfied and \eqref{A2TCs-diff} reduces to the flux equality, which is equivalent to
\begin{equation} \label{A2IF-diff-M1}
\begin{array}{c} \ds\int_{t^{n}}^{t^{n+1}}\ds\int_{E} \iF_{1}\iD_{1}(\lambda_{a}, \lambda, f,c_{0}) - \iF_{2}\iD_{2}(\lambda_{a}, \lambda, f,c_{0}) = 0,\vspace{0.15cm} \\
\hspace{4cm}\forall  E \in \iG_{h}, \; \; \forall  n=0, \hdots, N-1.
\end{array} 
\end{equation}
Note that the composite operator $ \iF_{i} \iD_{i}, \; i=1,2, $ is a Steklov-Poincar\'e (Dirichlet-to-Neumann) type operator. Equation \eqref{A2IF-diff-M1} together with \eqref{A2IF-adv-M1}  forms an interface problem, equivalent to
Problem~\ref{algo:splittingmulti}:
\begin{equation} \label{A2IFfull-M1} 
\begin{array}{ll}
\text{Find $ \left ( \lambda_{a}, \lambda\right ) \in (\lh)^{N \times L} \times (\lh)^{N} $ such that} \vspace{0.15cm} \\
\ds\int_{t^{n,l}}^{t^{n,l+1}} \ds\int_{E}\lambda_{a} -  \iH_{1}\iD_{1}(\lambda_{a}, \lambda, f,c_{0}) -  \iH_{2}\iD_{2}(\lambda_{a}, \lambda, f,c_{0}) &= 0, \vspace{0.15cm}\\
\ds\int_{t^{n}}^{t^{n+1}} \ds\int_{E} \iF_{1}\iD_{1}(\lambda_{a}, \lambda, f,c_{0}) + \iF_{2}\iD_{2}(\lambda_{a}, \lambda, f,c_{0}) & = 0,  \vspace{0.15cm} \\
& \hspace{-6cm}\forall E \in \iG_{h}, \; \; \forall  n=0, \hdots, N-1, \;\; \forall l=0, \hdots, L-1,
\end{array}
\end{equation}
or equivalently
\begin{equation} \label{A2IF-M1}
\begin{array}{c}
\text{Find $ \left ( \lambda_{a}, \lambda\right ) \in (\lh)^{N \times L} \times (\lh)^{N} $ such that} \vspace{0.15cm} \\
\iS \left (\begin{array}{c} \lambda_{a} \\
\lambda
\end{array}\right ) = \left (\begin{array}{c} \check{\chi } \\
\chi
\end{array}\right ),
\end{array}
\end{equation}
where 
\begin{equation*}
\iS \left (\begin{array}{c} \lambda_{a} \\
\lambda
\end{array}\right )  =  \left (\begin{array}{l} \ds\int_{t^{n,l}}^{t^{n,l+1}} \ds\int_{E} \lambda_{a} -  \sum_{i=1}^{2}\iH_{i}\iD_{i}(\lambda_{a}, \lambda, 0,0) \\[4mm]
\ds\int_{t^{n}}^{t^{n+1}} \ds\int_{E} -\sum_{i=1}^{2}\iF_{i}\iD_{i}(\lambda_{a}, \lambda, 0,0)
\end{array}\right )_{E \in \iG_{h}, \; n=0, \hdots, N-1 , \; l=0, \hdots, L-1}
\end{equation*}
and
\begin{equation*}
\left (\begin{array}{c} \check{\chi } \\
\chi
\end{array}\right ) = \left (\begin{array}{c} \ds\int_{t^{n,l}}^{t^{n,l+1}} \ds\int_{E} \sum_{i=1}^{2} \iH_{i}\iD_{i}(0, 0, f,c_{0}) \\[4mm]
\ds\int_{t^{n}}^{t^{n+1}} \ds\int_{E} \sum_{i=1}^{2}\iF_{i}\iD_{i}(0, 0, f,c_{0})
\end{array}\right )_{E \in \iG_{h}, \; n=0, \hdots, N-1 , \; l=0, \hdots, L-1}.
\end{equation*}

\noindent\\[-3mm]
\indent
System \eqref{A2IF-M1} can be solved iteratively by using a Krylov method (e.g. GMRES): the right hand side is computed only once by solving Problem~\ref{algo:splittingD} in each subdomain with $ \lambda_{a} = 0$ and $\lambda = 0 $; then for a pair of vectors $ (\eta_{a}, \eta) $ given in $ (\lh)^{N \times L} \times (\lh)^{N} $, the matrix vector product is obtained, at each Krylov iteration, by solving subdomain Problem~\ref{algo:splittingD} with $ \lambda_{a} = \eta_{a} $, $ \lambda = \eta $ and with $ f =0 $ and $ c_{0} = 0 $, and extracting the correct traces on the interface. 

  \begin{remark}
    Note that the method in~\cite{gander2014dirichlet} is related to the one introduced in~\cite{PhuongSINUM} and used here. It corresponds to a Richardson iteration applied to the above interface problem (with diffusion only) or that of Method 1 in~\cite{PhuongSINUM}. In the present paper as well as in~\cite{PhuongSINUM}, we use GMRES instead of Richardson, as it gives faster convergence for stationary problems. However, the situation is different for time dependent problems as pointed out and analyzed in~\cite{Nevanlinna} and one should use the convolution Krylov subspace methods for dynamical systems to accelerate the convergence to the same degree as in the case of stationary problems. An additional advantage in using GMRES is that it does not require a relaxation parameter (for which no optimal choice is known in the case of advection-diffusion problems) as does the method of~\cite{gander2014dirichlet}.
 \end{remark}

Following the same idea as in~\cite{PhuongSINUM} we apply a generalized Neumann-Neumann preconditioner. With this aim, we define the solution operator $ \iN_{i}, \; i=1,2: $
\begin{equation*} 
\begin{array}{rl} \iN_{i}: \left (\lh\right) ^{N\times L}  \times \left (\lh\right) ^{N}  &\rightarrow  \left (N_{h,i}\right) ^{N\times L}  \times \left (M_{h,i}\right) ^{N} \times \left (\Sigma_{h,i}\right) ^{N} \vspace{0.1cm}\\
\left (\lambda_{a}, \varphi \right ) & \mapsto \left (\hat{c}_{h,i}^{\Dt, \Dt_{a}}, c_{h,i}^{\Dt}, \br_{h,i}^{\Dt}\right ), 
\end{array}
\end{equation*} 
where 
$$ \hat{c}_{h,i}^{\Dt, \Dt_{a}} = \left (\hat{c}_{h,i}^{n,l}\right )_{n=0, \hdots, N-1, \; l=0, \hdots, L-1} \; \; \text{and} \; \;  \left (c_{h,i}^{\Dt}, \br_{h,i}^{\Dt}\right ) = \left (c_{h,i}^{n}, \br_{h,i}^{n}\right )_{n=1, \hdots, N} $$
are the solution of the subdomain problem that consists of solving, for~$ n~=~0, \hdots, N-~1 $,
\begin{itemize}
\item the advection equation: for $ l=0,\hdots, L-1 $,
\begin{equation*} \label{A2time-advection-PrecondM1}
			\begin{array}{rcll} \ds\int_{K} \phi_{i} \frac{c_{h,i}^{n,l+1} - c_{h,i}^{n,l} }{\Dt_{a}} + \sum_{E \subset \partial K} \ds\int_{E} \hat{c}_{h, i}^{n,l} (\bu \cdot \bn_{K}) &=&0, & \forall K \in \iK_{h,i}, \\
			\hat{c}_{h,i}^{n,l} & = &\iU_{h,i} \left (c_{h,i}^{n,l}, \lha^{n,l}\right ), 
			\end{array}
			\end{equation*}
			with $ c_{h,i}^{n,0} := c_{h,i}^{n} $ where $ c_{h,i}^{0}:= 0, $ 
\item and the diffusion equation
			\begin{equation*} \label{A2time-diffusion-PrecondM1}
			\hspace{-1cm}\begin{array}{rll} \ds\int_{K} \phi_{i} \frac{c_{h,i}^{n+1}-c_{h,i}^{n, L} }{\Dt} + \ds\int_{K} \Div \br_{h,i}^{n+1} & = 0, & \forall K \in \iK_{h,i} , \vspace{0.15cm}\\
			\ds\int_{\Omega_{i}} \bD_{i}^{-1} \br_{h,i}^{n+1} \cdot \bv -\ds\int_{\Omega_{i}}  c_{h,i}^{n+1} \Div \bv &=0,  & \forall \bv \in \Sigma_{h,i}^{0} ,\vspace{0.15cm}\\
			\ds\int_{E} \br_{h,i}^{n+1} \cdot \bn_{i} & = \ds\int_{E} \varphi^{n+1}, & \forall E \in \iG_{h}, \\
			\end{array}
			\end{equation*}
                        
	\end{itemize}
where $ \Sigma_{h,i}^{0}:= \left \{ \bv \in \Sigma_{h,i}: \bv \cdot \bn_{\mid E} = 0, \; \forall E \in \iG_{h} \right \}$
                         is introduced to treat the Neumann boundary condition on the interface.
%
%
                         
In order to define a (pseudo-)inverse operator of $ \iF_{i} \iD_{i}, \; i=1,2, $ we need to introduce the trace operator
\begin{equation*} 
\begin{array}{rl} \text{Tr}_{i}: \left (N_{h,i}\right) ^{N\times L}  \times \left (M_{h,i}\right) ^{N} \times \left (\Sigma_{h,i}\right) ^{N} & \rightarrow \left (\lh\right) ^{N}    \vspace{0.1cm}\\
\left (\hat{c}_{h,i}^{\Dt, \Dt_{a}}, c_{h,i}^{\Dt}, \br_{h,i}^{\Dt}\right ),  & \mapsto \underline{\lambda}
\end{array}
\end{equation*}
where $ \underline{\lambda}=\left (\underline{\lambda}_{h}^{n}\right )_{n=1, \hdots, N} $ stands for the trace of the concentration on the interface and is defined by
$$ \ds\int_{E} \underline{\lambda}_{h}^{n} (\bv_{E} \cdot \bn_{i}) = \ds\int_{\Omega_{i}} \bD_{i}^{-1} \br_{h,i}^{n} \cdot \bv_{E} -\ds\int_{\Omega_{i}}  c_{h,i}^{n} \Div \bv_{E}, \; \; \forall E \in \iG_{h}, \; n=1, \hdots, N,
$$
for $ \bv_{E} \in \Sigma_{h,i} $ such that $ \left (\bv_{E} \right )_{\mid {K}} = 0 $ for all $ K \in \iK_{h,i} $ that do not share the edge~$E$.

With these operators in place, the action of the generalized Neumann-Neumann preconditioner for \eqref{A2IF-M1}
on $ \left ( \mu_{a}, \mu\right ) \in (\lh)^{N \times L} \times (\lh)^{N},$ is defined by
\begin{equation*}
  \left(\begin{array}{c}
    \mu_{a}\\
    \mu
  \end{array}
  \right)
    \mapsto
    \left(\begin{array}{c}
\ds\int_{t^{n,l}}^{t^{n,l+1}} \ds\int_{E} \mu_{a} -  \sum_{i=1}^{2} \iH_{i}\iN_{i} (\mu_{a}, \mu)  \\
\ds\int_{t^{n}}^{t^{n+1}} \ds\int_{E} \sum_{i=1}^{2}\sigma_{i} \text{Tr}_{i}\iN_{i}(\mu_{a}, \mu)
\end{array}
\right)_{E \in \iG_{h}, \; \; \forall n=0, \hdots, N-1, \; \; \forall  l=0, \hdots, L-1} 
\end{equation*}
Here the composite operator $ \text{Tr}_{i}\iN_{i}, \; i=1,2, $ is a Neuman-to-Dirichlet type operator (which is the inverse operator of $ \iF_{i} \iD_{i}, \; i=1,2$) and $ \sigma_{i}: \lh \rightarrow [0,1], i=1,2, $ are weights such that $ \sigma_{1}+\sigma_{2} = 1 $. As in the case of pure diffusion problems (see~\cite{PhuongThesis,PhuongSINUM}), if $ \bD_{i} = d_{i} \mathbf{I}, i=1,2, $ then
$ \displaystyle\sigma_{i} := \frac{d_{i}}{d_{1} +d_{2}}. 
$
%
%
%
%
%
%
\subsection{Global-in-time optimized Schwarz (GTO-Schwarz) : an extension of the Optimized Schwarz Waveform Relaxation approach}
\label{A2Sec:M2IF}
As for the GTP-Schur method, we first define several operators needed to define the interface problem for this method. Let $ \iR_{i}:=\iLRi, \; i=1,2, $ be the solution operator
 which depends on the Robin parameter $ \prm_{i,j}, \; i=1,2, \; j=3-i: $
\begin{equation*} 
\begin{array}{rl} \iR_{i}: \left (\lh\right) ^{N\times L}  \times  \left (\lh\right) ^{N} \times & \hspace{-0.3cm}  L^{2} (0,T;L^{2}(\Omega_{i})) \times H_{*}^{1}(\Omega_{i}) \\
& \rightarrow  \left (\lh\right) ^{N} \times  \left (N_{h,i}\right) ^{N\times L}  \times \left (M_{h,i}\right) ^{N} \times \left (\Sigma_{h,i}\right) ^{N} \vspace{0.1cm}\\
\left (\lambda_{a}, \xi, f, c_{0}\right )  & \mapsto \left (\xi, \hat{c}_{h,i}^{\Dt, \Dt_{a}}, c_{h,i}^{\Dt}, \br_{h,i}^{\Dt}\right ):=\iLRi\left (\lambda_{a}, \xi, f, c_{0}\right ) , 
\end{array}
\end{equation*}
where $\iLRi$ is defined in~\eqref{eq:LR} and
\begin{itemize}
\item $ \lambda_{a} = \left (\lambda_{h,a}\right )_{n=0, \hdots, N-1, \; l=0, \hdots, L-1} $ represents Dirichlet boundary data on the interface for the advection equation (just as for the GTP-Schur method). 
\item $ \xi = \left (\xi_{h}^{n}\right )_{n=1, \hdots, N} $ represents the Robin boundary data (instead of Dirichlet data as for the GTP-Schur method) on the interface for the diffusion equation. Here we include $ \xi $ in the output of $ \iD_{i} $ as in the pure diffusion case (see~\cite{PhuongThesis,PhuongSINUM}) in order to compute Robin data transmitted to the neighboring subdomain.
\end{itemize} 
As stated in Remark~\ref{A2rmrkDiffStepAdvstep2M}, the advection step for the GTP-Schur method and the GTO-Schwarz method are the same. So we define the projection operator $ \widetilde{\iH}_{i}, \; i=1,2, $ similar to the operator $ \iH_{i} $ in the GTP-Schur method, but that takes as input the second component of the output of $ \iR_{i} $, instead
of its first component:
\begin{equation*} 
\hspace{-0.3cm}\begin{array}{rl} \widetilde{\iH}_{i}: \left (\lh\right) ^{N} \times \left (N_{h,i}\right) ^{N\times L}  \times \left (M_{h,i}\right) ^{N} \times \left (\Sigma_{h,i}\right) ^{N} & \rightarrow  \left (\lh\right) ^{N\times L}    \vspace{0.15cm}\\
\left (\xi, \hat{c}_{h,i}^{\Dt, \Dt_{a}}, c_{h,i}^{\Dt}, \br_{h,i}^{\Dt}\right )  & \mapsto \left \{ \begin{array}{cl} 
0, & \forall E \in \iG_{h,i}^{\text{in}}, \vspace{0.15cm}\\
\left (\hat{c}_{h,i}^{\Dt, \Dt_{a}}\right )_{\mid {E}} , & \forall E \in \iG_{h,i}^{\text{out}}. 
\end{array} \right .
\end{array} \vspace{0.1cm}
\end{equation*}

Next, for the Robin transmission conditions \eqref{A2RobinTCs-diff} of the diffusion equation, we need the following interface operators defined for $ i=1,2,$ with $ j=3-i, $
\begin{equation*} 
\begin{array}{l} \iB_{i}: \left (\lh\right) ^{N} \times  \left (N_{h,i}\right) ^{N\times L}  \times \left (M_{h,i}\right) ^{N} \times \left (\Sigma_{h,i}\right) ^{N}  \rightarrow \left (\lh\right) ^{N}    \vspace{0.1cm}\\
\left (\xi, \hat{c}_{h,i}^{\Dt, \Dt_{a}}, c_{h,i}^{\Dt}, \br_{h,i}^{\Dt}\right ) \mapsto \left (\br_{h,i}^{\Dt} \cdot \bn_{i} +\frac{\prm_{j,i}}{\prm_{i,j}} (\xi +\br_{h,i}^{\Dt} \cdot \bn_{i})\right )_{\mid E}, \; \forall E \in \iG_{h}.
\end{array}
\end{equation*}  
The transmission condition \eqref{A2TCs-adv} for the advection part leads to 
\begin{equation} \label{A2IF-adv-M2}
\begin{array}{c}
\ds\int_{t^{n,l}}^{t^{n,l+1}} \ds\int_{E} \lambda_{a} - \widetilde{\iH}_{1}\iR_{1}(\lambda_{a}, \xi_{1}, f,c_{0}) - \widetilde{\iH}_{2}\iR_{2}(\lambda_{a}, \xi_{2}, f,c_{0}) = 0, \vspace{0.15cm}\\
\hspace{2cm} \forall E \in \iG_{h}, \; \; \forall n=0, \hdots, N-1, \; \forall l=0, \hdots, L-1.
\end{array}
\end{equation}

Exploiting the boundary condition~\eqref{eq:RobinContinu}, we see that \eqref{A2RobinTCs-diff} is equivalent to
\begin{equation} \label{A2IF-diff-M2}
\begin{array}{l} \ds\int_{t^{n}}^{t^{n+1}} \ds\int_{E} \xi_{1} - \iB_{2} \iR_{2} (\lambda_{a}, \xi_{2}, f,c_{0})  = 0, \vspace{0.15cm}\\
\ds\int_{t^{n}}^{t^{n+1}} \ds\int_{E} \xi_{2} - \iB_{1} \iR_{1} (\lambda_{a}, \xi_{1}, f,c_{0})  = 0, 
\end{array} \; \forall E \in \iG_{h}, \; \; \forall n=0, \hdots, N-1. 
\end{equation}
Note that the composite operator $ \iB_{i} \iR_{i}, \; i=1,2, $ is a discrete Robin-to-Robin type operator. Equation \eqref{A2IF-diff-M2} together with \eqref{A2IF-adv-M2}  forms an interface problem, equivalent to
Problem~\ref{algo:splittingmulti} were we have replaced \eqref{A2TCs-diff} by \eqref{A2RobinTCs-diff}, as follows
\begin{equation} \label{A2IFfull-M2} 
\begin{array}{ll}
\text{Find $ \left (\lambda_{a}, \xi_{1}, \xi_{2}\right ) \in (\lh)^{N \times L} \times (\lh)^{N} \times (\lh)^{N} $ such that} \vspace{0.15cm}\\
\ds\int_{t^{n,l}}^{t^{n,l+1}}  \ds\int_{E}\lambda_{a} - \widetilde{\iH}_{1}\iR_{1}(\lambda_{a}, \xi_{1}, f,c_{0}) - \widetilde{\iH}_{2}\iR_{2}(\lambda_{a}, \xi_{2}, f,c_{0}) &= 0,\vspace{0.15cm}\\
\ds\int_{t^{n}}^{t^{n+1}} \ds\int_{E} \xi_{1} - \iB_{1} \iR_{2} (\lambda_{a}, \xi_{2}, f,c_{0})  & = 0, \vspace{0.15cm}\\
\ds\int_{t^{n}}^{t^{n+1}} \ds\int_{E} \xi_{2} - \iB_{2} \iR_{1} (\lambda_{a}, \xi_{1}, f,c_{0})  & = 0,\vspace{0.15cm}\\
& \hspace{-8cm} \forall E \in \iG_{h}, \; \; \forall n=0, \hdots, N-1, \; \forall l=0, \hdots, L-1,
\end{array}
\end{equation}
or equivalently,
\begin{equation} \label{A2IF-M2}
\iS_{R} \left (\begin{array}{c} \check{ \lambda} \\
\xi_{1} \\ \xi_{2} 
\end{array} \right ) = \chi_{R}, \\
\end{equation}
where 
\begin{equation*}
\iS_{R}  \left (\begin{array}{c} \check{ \lambda} \\
\xi_{1} \\ \xi_{2} 
\end{array} \right ) =  \left (\begin{array}{l} \ds\int_{t^{n,l}}^{t^{n,l+1}}  \ds\int_{E} \lambda_{a} - \sum_{i=1}^{2} \widetilde{\iH}_{i}\iR_{i}(\lambda_{a}, \xi_{i}, 0,0) \vspace{0.15cm}\\
\ds\int_{t^{n}}^{t^{n+1}}  \ds\int_{E} \xi_{1} - \iB_{1} \iR_{2} (\lambda_{a}, \xi_{2}, 0,0) \vspace{0.15cm}\\
\ds\int_{t^{n}}^{t^{n+1}}  \ds\int_{E} \xi_{2} -  \iB_{2} \iR_{1} (\lambda_{a}, \xi_{1}, 0, 0)
\end{array} \right )_{E \in \iG_{h}, \; n=0, \hdots, N-1 , \; l=0, \hdots, L-1}
\end{equation*}
and
\begin{equation*}
\chi_{R} = \left (\begin{array}{l} \ds\int_{t^{n,l}}^{t^{n,l+1}}  \ds\int_{E} \sum_{i=1}^{2} \widetilde{\iH}_{i}\iR_{i}(0, 0, f,c_{0}) \vspace{0.15cm}\\
\ds\int_{t^{n}}^{t^{n+1}}  \ds\int_{E} \iB_{1} \iR_{2} (0, 0, f,c_{0}) \vspace{0.15cm}\\ 
\ds\int_{t^{n}}^{t^{n+1}}  \ds\int_{E} \iB_{2} \iR_{1} (0, 0, f,c_{0})
\end{array} \right )_{E \in \iG_{h}, \; n=0, \hdots, N-1 , \; l=0, \hdots, L-1}.
\end{equation*}

System~\eqref{A2IF-M2} can be solved iteratively using Jacobi iterations (which corresponds to the discrete "splitting" OSWR algorithm) or a Krylov method such as GMRES: the right hand side is computed by solving
Problem~\ref{algo:splittingR} in each subdomain with $ \lambda_{a}~=~0$ and $\xi = 0 $; then for a given vector $ (\eta_{a}, \xi_{1}, \xi_{2}) $ in $ (\lh)^{N \times L} \times (\lh)^{N}  \times (\lh)^{N} $, the matrix vector product is obtained (at each iteration) by solving Problem~\ref{algo:splittingR} in $ \Omega_{i} \times (0,T), \; i=1,2, $ with $ \lambda_{a} = \eta_{a} $, $ \xi = \xi_{i}, $ and with $ f =0 $ and $ c_{0} = 0 $.
\begin{remark}
Due to the use of the splitting method, we have formulated a generalization of the OSWR method in which the Robin parameters only act on the diffusion equation as in the case of pure diffusion problems. The advection term is now like a source term for the diffusion equation. Thus the optimized Robin parameters $ \prm_{i,j},$ $ i~=~1,2,$ $j~=~3~-~i, $ are calculated in the same way as for the pure diffusion case (see~\cite{PhuongSINUM}). Consequently, the advection coefficient is not taken into account in the computation of the optimized parameters. This may be an advantage of using operator splitting because we don't need to handle variable coefficients due to the velocity field, as one does for a fully implicit scheme (see, e.g, \cite{OSWRDG, OSWRDG2}). In Section~\ref{A2Sec:Num}, we study numerically the impact of the optimized parameters on the convergence behavior, especially for advection-dominated problems. 
\end{remark}
As the interface problem derived above for each method is global in time, one may use different time steps for different subdomains as in the case of pure diffusion problems (see~\cite{PhuongSINUM}). In the next section, we describe how we enforce the transmission conditions over such nonconforming time grids. 
%
%
%
\section{Nonconforming time discretizations}
\label{Sec:AdvTime}
\begin{figure}[h]
\vspace{2cm}
\hspace{2.cm}
\begin{minipage}[c]{0.5 \linewidth}
\setlength{\unitlength}{1pt} 
\begin{picture}(140,70)(0,0)
\thicklines
\put(0,3){\includegraphics[scale=0.7]{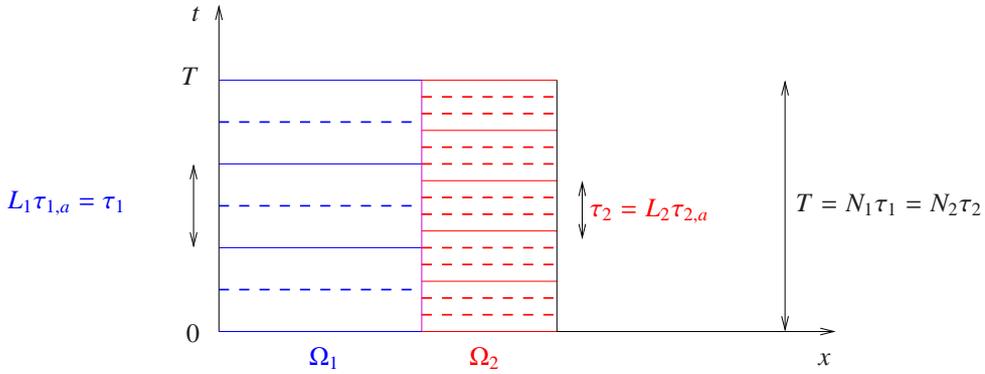} \\}
\put(-1,1){$ 0 $}
\put(-3,98){$ T $}
\put(45,-8){\textcolor{blue}{$ \Omega_{1} $}}
\put(105,-8){\textcolor{red}{$ \Omega_{2}$}}
\put(-68,51){\textcolor{blue}{$L_{1} \Dt_{1,a} = \Dt_{1} $}}
\put(150,48){\textcolor{red}{$ \Dt_{2} =L_{2} \Dt_{2,a} $}}
\put(227,50){{$ T = N_{1} \Dt_{1} = N_{2} \Dt_{2}$}}
\put(235,-8){$ x $}
\put(1,122){$ t $}
\end{picture}
\end{minipage} \vspace{0.5cm}
\caption{Nonconforming time grids in the subdomains.}
\label{A2Fig:AdvTime} 
\end{figure}

Let $ \iT_{1} $ and $ \iT_{2} $
be two different uniform partitions of the time interval $ (0,T) $ into $ N_{1} $ and $ N_{2} $ sub-intervals respectively with lengths $ \Dt_{1} $ and $ \Dt_{2} $, respectively (see Figure~\ref{A2Fig:AdvTime}). The sub-time step for the advection in each subdomain is defined by 
$ \Dt_{i} = L_{i} \Dt_{i,a}, \; i=1,2, 
$
and we denote by $ \iT_{i}^{a}, \; i=1,2, $ the corresponding partition in time for the advection. We denote by $ P_{0}(\iT_{i}, \lh) $ the space of piecewise constant functions in time on grid $ \iT_{i} $ with values in $ \lh $. Then define $ \Pi_{ij} $ the average-valued projection from $ P_{0}(\iT_{j}, \lh) $ to $ P_{0}(\iT_{i}, \lh) $ (see~\cite{PhuongSINUM})  and $ \Pi_{ij}^{a} $ from $ P_{0}(\iT_{j}^{a}, \lh) $ to $ P_{0}(\iT_{i}^{a}, \lh) $. 

As pointed out earlier, due to the use of the splitting method the interface problem for either GTP Schur or GTO Schwarz consists of an equation imposing the transmission condition for the advection problem and one or more imposing the transmission conditions for the diffusion problem. The latter can be enforced in time in a way similar to that of the pure diffusion problem (see~\cite{PhuongSINUM}). For the advection transmission condition, as there is only one unknown $ \lambda_{a} $ on the interface, one needs to choose $ \lambda_{a} $ to be piecewise constant in time on either the grid $ \iT_{1}^{a} $ or $ \iT_{2}^{a} $. This is the case for both methods. \vspace{15pt}

\textbf{For GTP Schur} \hspace{0.2cm}
We choose $ \lambda_{a} $ and $ \lambda $ to be piecewise constant in time on the advection and diffusion time grids respectively. For instance, let $ \lambda_{a} \in  P_{0}(\iT_{2}^{a}, \lh) $ and $ \lambda \in  P_{0}(\iT_{2}, \lh)$. Then the interface problem \eqref{A2IFfull-M1} is rewritten as
\begin{equation} \label{A2IF-M1-nonTime}
\begin{array}{ll}
\text{Find $ \left (\lambda_{a}, \lambda \right ) \in (\lh)^{N_{2} \times L_{2}} \times (\lh)^{N} $ such that} \vspace{0.15cm}\\
\ds\int_{t_{2}^{n,l}}^{t_{2}^{n,l+1}} \ds\int_{E}\lambda_{a} - \Pi_{21}^{a}\left ( \iH_{1}\iD_{1}(\Pi_{12}^{a}(\lambda_{a}), \Pi_{12}(\lambda), f,c_{0})\right ) -  \iH_{2}\iD_{2}(\lambda_{a}, \lambda, f,c_{0}) &= 0, \vspace{0.15cm}\\
\ds\int_{t_{2}^{n}}^{t_{2}^{n+1}} \ds\int_{E} \Pi_{21}\left (\iF_{1}\iD_{1}(\Pi_{12}^{a}(\lambda_{a}), \Pi_{12}(\lambda), f,c_{0})\right ) - \iF_{2}\iD_{2}(\lambda_{a}, \lambda, f,c_{0}) & = 0,
\end{array}
\end{equation}
for $ \; \forall E \in \iG_{h} $ and $ \; \forall n=0, \hdots, N_{2}-1, \; \; \forall l=0, \hdots, L_{2}-1. $ 
%
%
\vspace{15pt}

\textbf{For GTO Schwarz} \hspace{0.2cm}
We choose $ \lambda_{a} $ to be piecewise constant in time on one grid, for instance, $ \iT_{2}^{a} $. For the two Robin terms $ \xi_{1} $ and $ \xi_{2} $, we use the same technique as in~\cite{PhuongSINUM}. The interface problem \eqref{A2IFfull-M2} is then rewritten as
\begin{equation} \label{A2IF-M2-nonTime}
\begin{array}{ll}
\text{Find $ \left (\lambda_{a}, \xi_{1}, \xi_{2}\right ) \in (\lh)^{N_{2} \times L_{2}} \times (\lh)^{N_{1}} \times (\lh)^{N_{2}} $ such that} \vspace{0.15cm}\\
\ds\int_{t_{2}^{n,l}}^{t_{2}^{n,l+1}} \ds\int_{E}\lambda_{a} - \Pi_{21}^{a} \left (\widetilde{\iH}_{1}\iR_{1}(\Pi_{12}^{a} (\lambda_{a}), \xi_{1}, f,c_{0}) - \widetilde{\iH}_{2}\iR_{2}(\lambda_{a}, \xi_{2}, f,c_{0})\right ) &= 0,\vspace{0.15cm}\\
\ds\int_{t_{1}^{m}}^{t_{1}^{m+1}} \ds\int_{E} \xi_{1} - \Pi_{12}\left (\iB_{1} \iR_{2} (\lambda_{a}, \xi_{2}, f,c_{0})\right )  & = 0, \vspace{0.15cm}\\
\ds\int_{t_{2}^{n}}^{t_{2}^{n+1}} \ds\int_{E} \xi_{2} - \Pi_{21}\left (\iB_{2} \iR_{1} (\Pi_{12}^{a}(\lambda_{a}), \xi_{1}, f,c_{0})\right )  & = 0,
\end{array}
\end{equation}
for $ \; \forall E \in \iG_{h} $, $\; \; \forall m=0, \hdots, N_{1}-1, $ and $ \; \forall n=0, \hdots, N_{2}-1, \; \; \forall l=0, \hdots, L_{2}-1. $ 

For conforming time grids, the two schemes defined by performing GMRES on the two interface problems \eqref{A2IF-M1-nonTime} and \eqref{A2IF-M2-nonTime} respectively converge to the same monodomain solution, while for the nonconforming case, these two schemes yield different solutions at convergence due to the use of different projection operators (this is also the case for pure diffusion problems studied in~\cite{PhuongSINUM}). In the next section we will carry out numerical experiments to investigate and compare the errors in time of the two methods. 
%
%
%
%
%
%
%
%
%
\section{Numerical results}
\label{A2Sec:Num}
We present 2D numerical experiments 
to illustrate the performance of the two methods formulated in the previous sections.  We consider an isotropic diffusion matrix $ \bD_{i} = d_{i} \pmb{I} $, where $ \pmb{I} $ is the 2D identity matrix. 
In Subsection~\ref{A2subsec:TestPeclet}, a simple test case with two subdomains is studied. The coefficients are constant in the subdomains and can be continuous or discontinuous across the interface. We verify the convergence behavior of the two methods for different P\'eclet numbers.  
In Subsection~\ref{A2subsec:CEA} we consider a test case that is a prototype for a nuclear waste repository simulation, in which the subdomains involved have different length scales (from $ 1 $m to $ 100 $m) and different physical properties. The convergence of the two methods for a decomposition into $ 9 $ subdomains is studied, and we analyze numerically the error in time when nonconforming grids are used. Time windows are employed to approximate the solution over long time intervals. In Subsection~\ref{A2subsec:ANDRA}, a test case for the simulation of the transport around a surface nuclear waste storage site is considered. In this case, the geometry of the domain is quite complex and layers with highly different physical properties are present. The domain is decomposed into $ 6 $ subdomains and time windows are also used.

\begin{remark} \label{rmrk:SubSolve}
One iteration of the GTP-Schur method with the preconditioner
costs twice as much as one iteration of the GTO-Schwarz method (in terms of number of subdomain solves).
Thus to compare the convergence of the two methods with GMRES, in the sequel we show the error in the concentration $ c $ and the vector field $ \br $ versus the number of subdomain solves (instead of the number of iterations). 
\end{remark}

%
%
\subsection{Piecewise constant coefficients}
\label{A2subsec:TestPeclet}
The computational domain $ \Omega$ is the unit square,
and the final time is $ T=1 $.  We split $ \Omega $ into two nonoverlapping subdomains
$ \Omega_{1} = (0,0.5) \times (0,1) $ and $ \Omega_{2} = (0.5,1) \times (0,1) $. 
Homogeneous Dirichlet boundary conditions are imposed on $ \partial \Omega $, the initial condition is 
\begin{equation} \label{A2Test1:IC}
c_{0}(x,y)=xy(1-x)(1-y)\exp(-100((x-0.2)^2+(y-0.2)^2)),
\end{equation}  
and the source term is 
\begin{equation} \label{A2Test1:source}
f(x,y,t)~=~\exp(-100((x-0.2)^2+(y-0.2)^2)).
\end{equation}
The porosity is $ \phi _{1} = \phi_{2} =\phi = 1 $. 
The advection and diffusion coefficients, $ \bu_{i} $ and $ d_{i}, i=1,2, $ given in Table~\ref{A2Tab:Test1Discont}, are constant in each subdomain but discontinuous across the interface. The global P\'eclet number and the CFL condition in each subdomain are also shown. We use nonconforming time grids $ \Dt_{1} \neq \Dt_{2} $,  but equal advection and diffusion time steps, $ \Dt_{a,i} = \Dt_{i}, i=1,2.  $ In space, we use a uniform rectangular mesh with size $ \Delta x_{1} = \Delta x_{2} = \Delta x = 1/100 $. 
\begin{table}[h]
\centering
\begin{tabular}{|l|l|l|l|l|l|l|}
  \hline
  									& $ \bu_{i} $								&$  d_{i}$			& $ \text{Pe}_{G}$ 	&$ \Dt_{i}  $	 & $ dt_{\text{CFL}} $		& $ \Dt_{a,i} $ \\ \hline
   $ \Omega_{1} $        & $ \left (0.5, 1\right ) $ 			& $ 0.02 $  		& $ \approx 10$		& $ 1/100 $       	 & $ 1/100 $					& $ 1/100$ 				\\ \hline
   $ \Omega_{2} $        & $ \left (0.5, 0.1\right ) $        & $ 0.002 $ 		& $ \approx 100 $  	& $ 1/75 $ 	  		 & $ 1/50 $					& $ 1/75 $				\\ \hline
\end{tabular}
\caption{Data for the discontinuous test case.}
\label{A2Tab:Test1Discont} \vspace{-0.2cm}
\end{table}
%
%
%

\begin{figure}[h]
\begin{flushleft}
\begin{minipage}[c]{0.45 \linewidth}
\begin{center}
\includegraphics[scale=0.26]{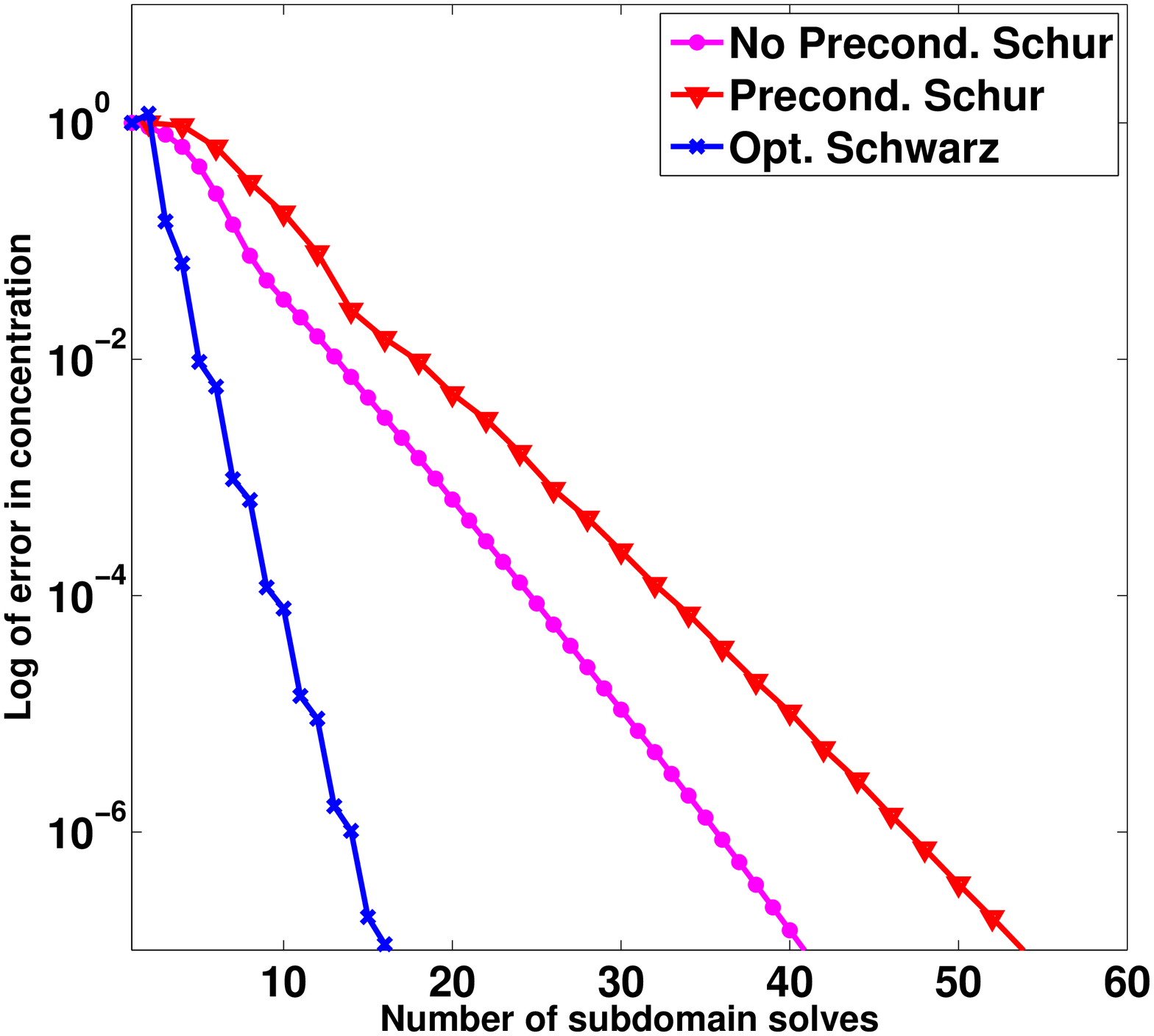}
\end{center}
\end{minipage} \hspace{10pt}
\begin{minipage}[c]{0.45 \linewidth}
\begin{center}
\includegraphics[scale=0.26]{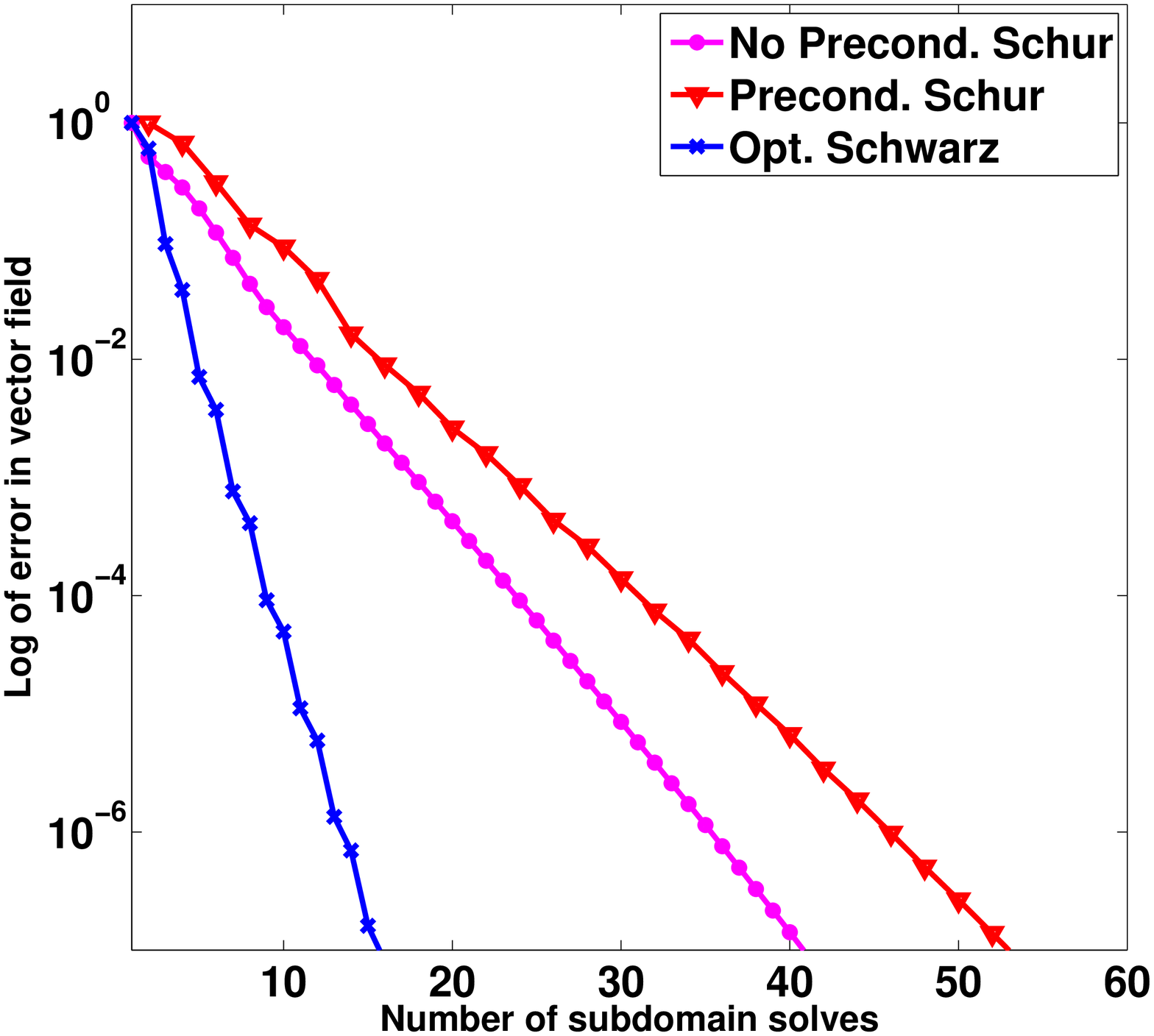}
\end{center}
\end{minipage}
\end{flushleft}
\caption{For discontinuous coefficients: Convergence curves for the different algorithms using GMRES: $ L^{2}-L^{2} $ error in $ c $ (left) and in $ \br $ (right). } 	
\label{A2Fig:Test1DiscontConvergence}  \vspace{-0.3cm}
\end{figure}

As in the continuous coefficient case, we analyze the convergence behavior of each method with~$ c_{0} = 0$ and $f = 0 $. Figure~\ref{A2Fig:Test1DiscontConvergence} shows the error (in logarithmic scale)
in the $ L^{2} (0,T; L^{2}(\Omega)) $-norm of the concentration $ c $ and of the vector field $ \br $, 
versus the number of subdomain solves using GMRES with a random initial guess. Again we see that for advection-dominated problems, the preconditioner for the GTP-Schur method does not work well and the GTO-Schwarz method converges much faster than the GTP-Schur method, by about a factor of $ 2.6 $ (with no preconditioner) and a factor of $ 3.3 $ (with the preconditioner) for both errors in $ c $ and~$ \br $. 

\begin{figure}[h]
\centering
\begin{minipage}[c]{0.5 \linewidth}
\includegraphics[scale=0.2]{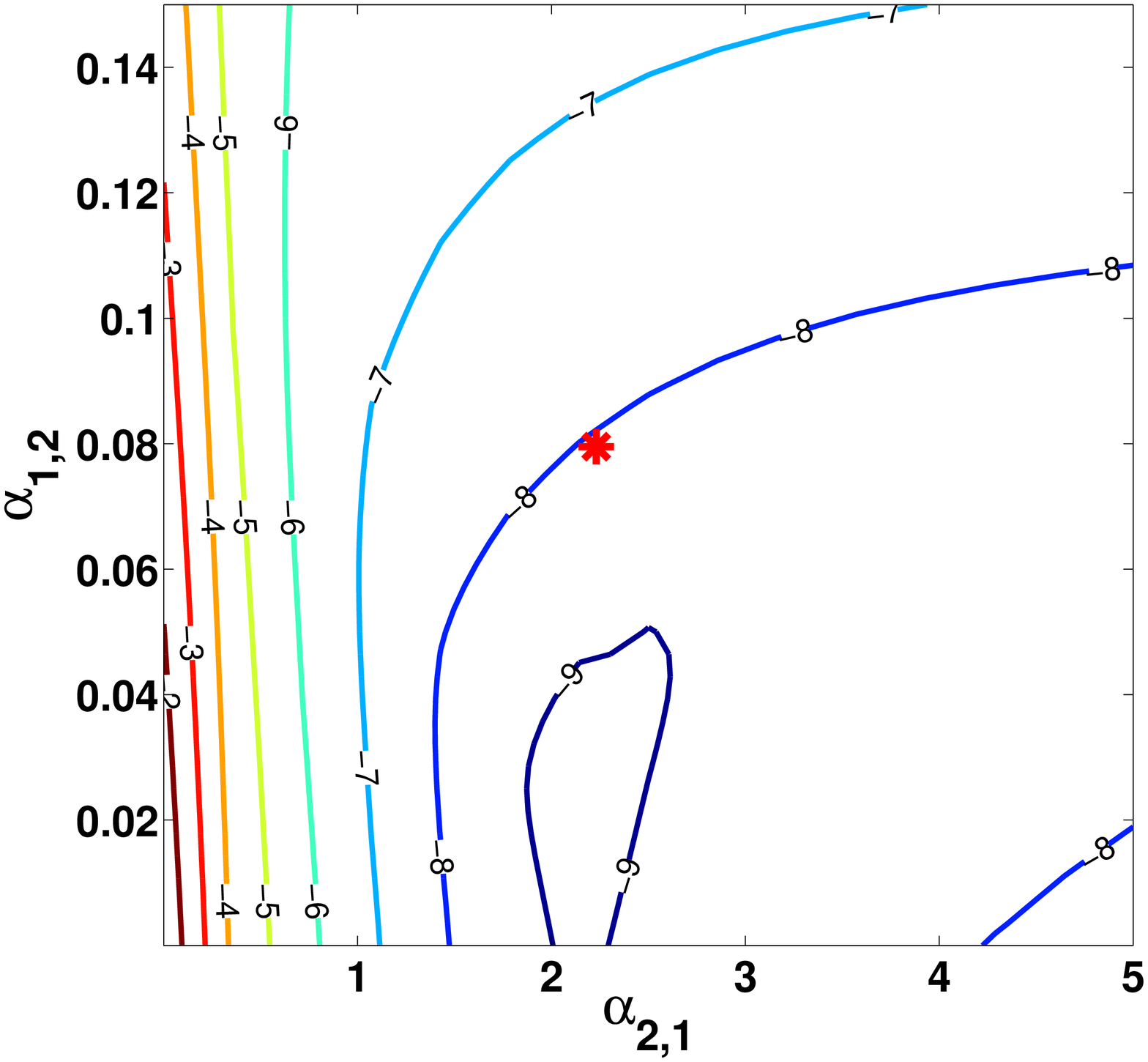}
\end{minipage} \vspace{-0.3cm}
\caption{For discontinuous coefficients: Level curves for the error in $ \br $ after $ 15 $ Jacobi iterations for various values of $ \alpha_{1,2} $ and $ \alpha_{2,1} $.} 	
\label{A2Fig:Test1DiscontRPara}  \vspace{-0.3cm}
\end{figure}

Figure~\ref{A2Fig:Test1DiscontRPara} shows the iso-lines of the 
error in the diffusive flux $ \br $ (in logarithmic scale) reched after $ 15 $ Jacobi iterations,
for various values of the parameters $ \alpha_{1,2} $ and $ \alpha_{2,1} $. We observe that, for discontinuous coefficients, the pair of optimized parameters (red star), computed as shown in~\cite{OSWR3sub}, is also located close to the optimal numerical values. 

%
%
%

In Table~\ref{A2Tab:Test1DisContAsymp}, the number of subdomain solves needed to reach a reduction of $ 10^{-6} $ in the error in the concentration $ c $ and in the vector field $ \br $ (in square brackets) when refining the mesh in space and in time with the ratio of $ \Delta x $ to $ \Dt $ constant is shown.
We see that for discontinuous coefficients and nonconforming time grids, the convergence of the GTP-Schur method with the Neumann-Neumann preconditioner is slightly dependent on the mesh size, while that with no preconditioner increases fairly rapidly with decreasing mesh size. For the GTO-Schwarz method, the convergence is almost independent of the mesh size and again, the use of GMRES, instead of Jacobi iterations, does not improve significantly the convergence speed. We see that the convergence of the GTO-Schwarz method is very fast. 
\begin{table}[h]
\centering
{\footnotesize
\begin{tabular}{|c|c|c|c|c|}
\hline 
\multirow{2}{*}{$ \Delta {x}$} 	& No Precond. 		& Precond.			&\multicolumn{2}{|c|}{Opt. Schwarz}  	\\ \cline{4-5}
								 				& Schur 				& Schur 			 	& GMRES				& Jacobi 							\\ \hline
$ 1/50 $									&	$ 25 \, [25]$		&	$ 40 \, [40]$		&	 $ 13 \, [13] $	&$ 14 \, [14]$								\\ \hline	
$ 1/100 $									&	$ 36 \, [36]$		&	$ 48 \, [46]$		&	 $ 15 \, [14] $	&$ 16 \, [15]$								\\ \hline				
$ 1/200 $									&	$ 52 \, [50]$		&	$ 54 \, [52]$		&	 $ 15 \, [15] $	&$ 17 \, [16]$										 \\ \hline	
$ 1/400 $									&	$ 76 \, [72]$		&	$ 58 \, [54]$		&	 $ 16 \, [15] $	&$ 18 \, [17]$									\\ \hline				
\end{tabular}}
\caption{For discontinuous coefficients: Number of subdomain solves required to reach a reduction of $ 10^{-6} $ in the error for the different algorithms , and for different values of the discretization parameters.}
\label{A2Tab:Test1DisContAsymp} 
\end{table}
%
%
%
%
%
%
%
%
%
%
%
	\subsection{A near-field simulation }
	\label{A2subsec:CEA}
\begin{figure}[h]
\begin{minipage}[c]{0.4 \linewidth}
\hspace{0.5cm}
\includegraphics[height=6.5cm]{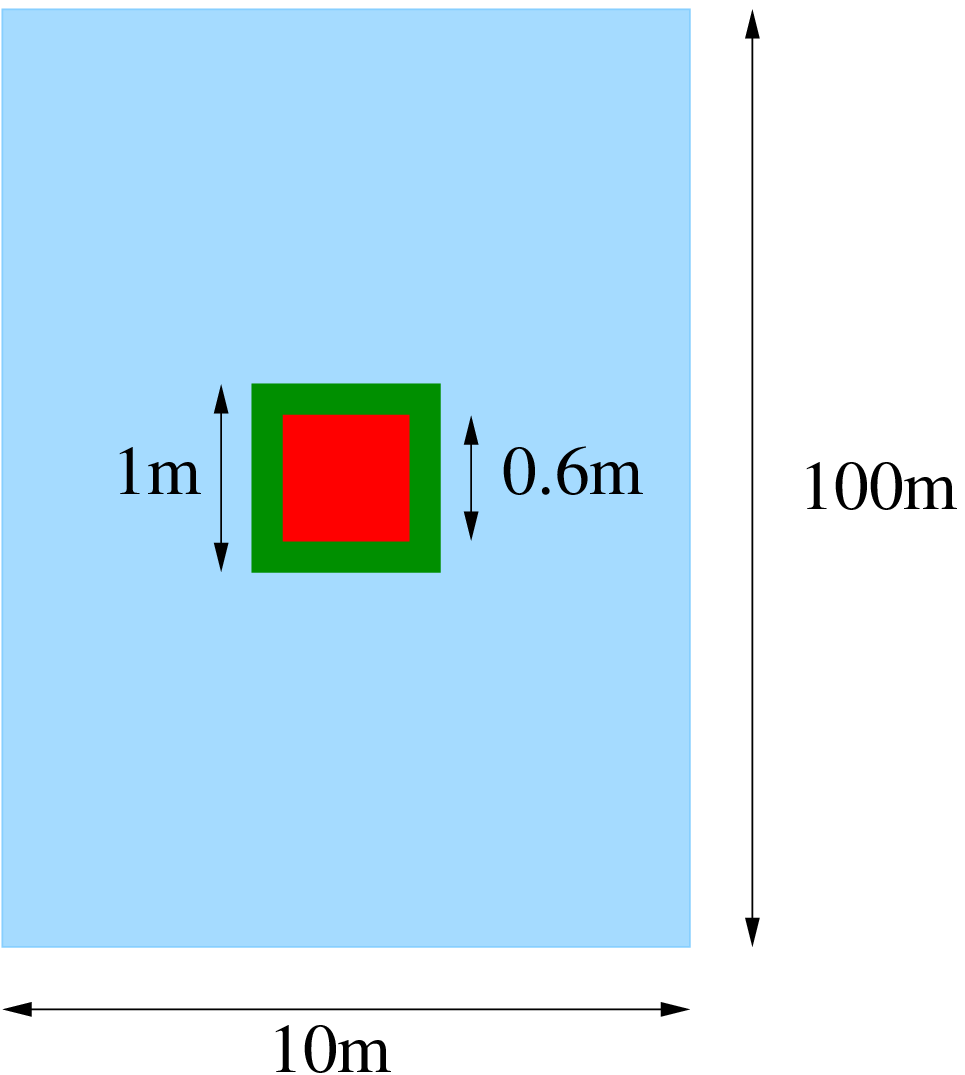} 
\end{minipage} \hspace{2.2cm}
\begin{minipage}[c]{0.3 \linewidth}
\vspace{5cm}
\setlength{\unitlength}{1pt} 
\begin{picture}(140,70)(0,0)
\thicklines
\put(0,22){\includegraphics[height=6.1cm]{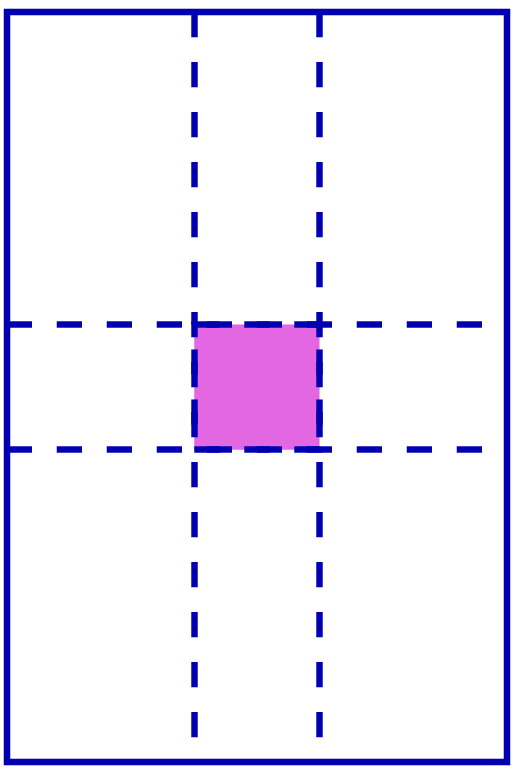} \\}
\put(17,60){$ \Omega_{1} $}
\put(52,60){$ \Omega_{2} $}
\put(91,60){$ \Omega_{3}$}
\put(20,106){$ \Omega_{4} $}
\put(52,106){$ \pmb{\Omega_{5}} $}
\put(95,106){$ \Omega_{6}$}
\put(17,170){$ \Omega_{7} $}
\put(52,170){$ \Omega_{8} $}
\put(91,170){$ \Omega_{9}$}
\end{picture}
\end{minipage} 
\caption{The domain of calculation and its decomposition. } 	
\label{A2Fig:CEAdomain} \vspace{-0.2cm}
\end{figure}

We consider a simplified test case~\cite{PaminaNearfield} for the simulation of the transport of a contaminant in a near field around a nuclear waste repository site.
The domain of calculation is a $ 10 $m by $ 100 $m rectangle and the repository is a centrally located unit square (see Figure \ref{A2Fig:CEAdomain} with a blow-up in the $ x- $ direction for visualization purpose). The repository consists of the EDZ (Excavation Damaged Zone) and the vitrified waste. The final time is $ T_{f} =  2 \times 10^{11} $s ($ \approx 20000 $ years). The coefficients for the simulation are given in Table \ref{A2Tab:CEAdata}. 
The advection field is governed by the (time-independent) Darcy flow equation together with the law of mass conservation
\begin{equation} \label{A2CEADarcyflow}
\begin{array}{rll} \bu & = - K \nabla p & \text{in} \; \Omega,\\
\Div \bu & = 0 & \text{in} \; \Omega,
\end{array} 
\end{equation}
where $p$ is the pressure field.
No flow boundary is imposed horizontally and a pressure gradient is imposed vertically with $ p=100 $ Pa  on bottom and $ p=0 $ on top. 

The source term is $ f=0 $ and an initial condition $ c_{0} $ is defined by
\begin{equation} \label{A2CEA:IC}
c_{0}= \left \{ \begin{array}{ll} 1, & \text{in the red box (containing the vitrified waste)}, \\
0, & \text{elsewhere}.
\end{array} \right. 
\end{equation}
Boundary conditions of the transport problem are homogeneous Dirichlet conditions on top and bottom, and homogeneous Neumann conditions on the left
and right hand sides. 
\begin{table}[h]
\centering
\begin{tabular}{|l|l|l|l|l|}
  \hline
	    		Material                 				& Permeability $ K $ (m/s)       	& Porosity  $ \phi $               & Diffusion coefficient $ d $ (m$^{2}$/s) 				 \\ \hline
	    		Host rock								& $ 10^{-13} $   				& $ 0.06 $                 & $ 6 \, 10^{-13} $    	 \\ \hline 
	    		EDZ										& $ 5 \, 10^{-11} $   		& $ 0.20 $                  & $2 \, 10^{-11} $     	  \\ \hline 	
	    		Vitrified waste						& $ 10^{-8} $   					& $ 0.10 $                  & $ 10^{-11} $     	  \\ \hline 	     		
\end{tabular}
\caption{Data for flow and transport problems.}
\label{A2Tab:CEAdata} \vspace{-0.2cm}
\end{table}

\begin{figure}[h]
\centering
\hspace{1.5cm}\begin{minipage}[c]{0.8 \linewidth}
\includegraphics[scale=0.2]{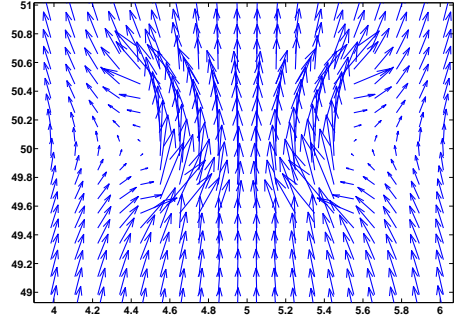} 
\end{minipage}  \vspace{-0.5cm}
\caption{Darcy flow. } 	
\label{A2Fig:Test3CEADarcyvelo} \vspace{-0.2cm}
\end{figure}

For the spatial discretization (for both the flow and transport equations), we use a non-uniform but conforming rectangular mesh with a finer discretization in the repository (a uniform mesh with $ 10 $ points in each direction) and a coarser discretization in the host rock (the mesh size progressively increases with distance from the repository by a factor of $ 1.05 $). The Darcy flow is approximated by using mixed finite elements. Figure~\ref{A2Fig:Test3CEADarcyvelo} shows a zoom of the velocity field around the repository. 
The maximum local P\'eclet number in this test case is $ 0.0513 $, thus it is a diffusion-dominated problem. The time step due to the CFL condition is large as the velocity field is very small (of order of $ 10^{-13} $ m/s), $ \Dt_{\text{CFL}} = 0.075 \; T_{f} $ in the repository and $ \Dt_{\text{CFL}}  = 0.125 \; T_{f}$ elsewhere. 
We decompose $ \Omega $ into $ 9 $ subdomains as depicted in Figure~\ref{A2Fig:CEAdomain} with $ \Omega_{5} $ representing the repository.
For the time discretization, we use nonconforming time grids (with a finer time step in the repository) and equal diffusion and advection time steps 
$\Dt_{i} = \Dt_{a,i} , \forall i. $

We observe that the longer the time interval the slower the convergence. In addition, for a fixed time step $ \Dt $, it is more costly to approximate the solution for a longer time interval than for a shorter time interval. Thus we use time windows for this test case. We divide $ (0,T_{f}) $ into $ 200 $ time windows with size $ T = 10^{9} $s. We will first analyze the convergence behavior as well as the accuracy in time of the multidomain solution with nonconforming grids for the first time window, $ (0,T) $. The time steps are $ \Dt_{5} = \Dt_{a,5} = T/500, $ and $\Dt_{i} = \Dt_{a,i} = T/100 $, $ i \neq 5 $. 
%
%

\begin{figure}[h]
\centering
\hspace{-7mm}
\begin{minipage}{0.45 \linewidth}
\includegraphics[scale=0.24]{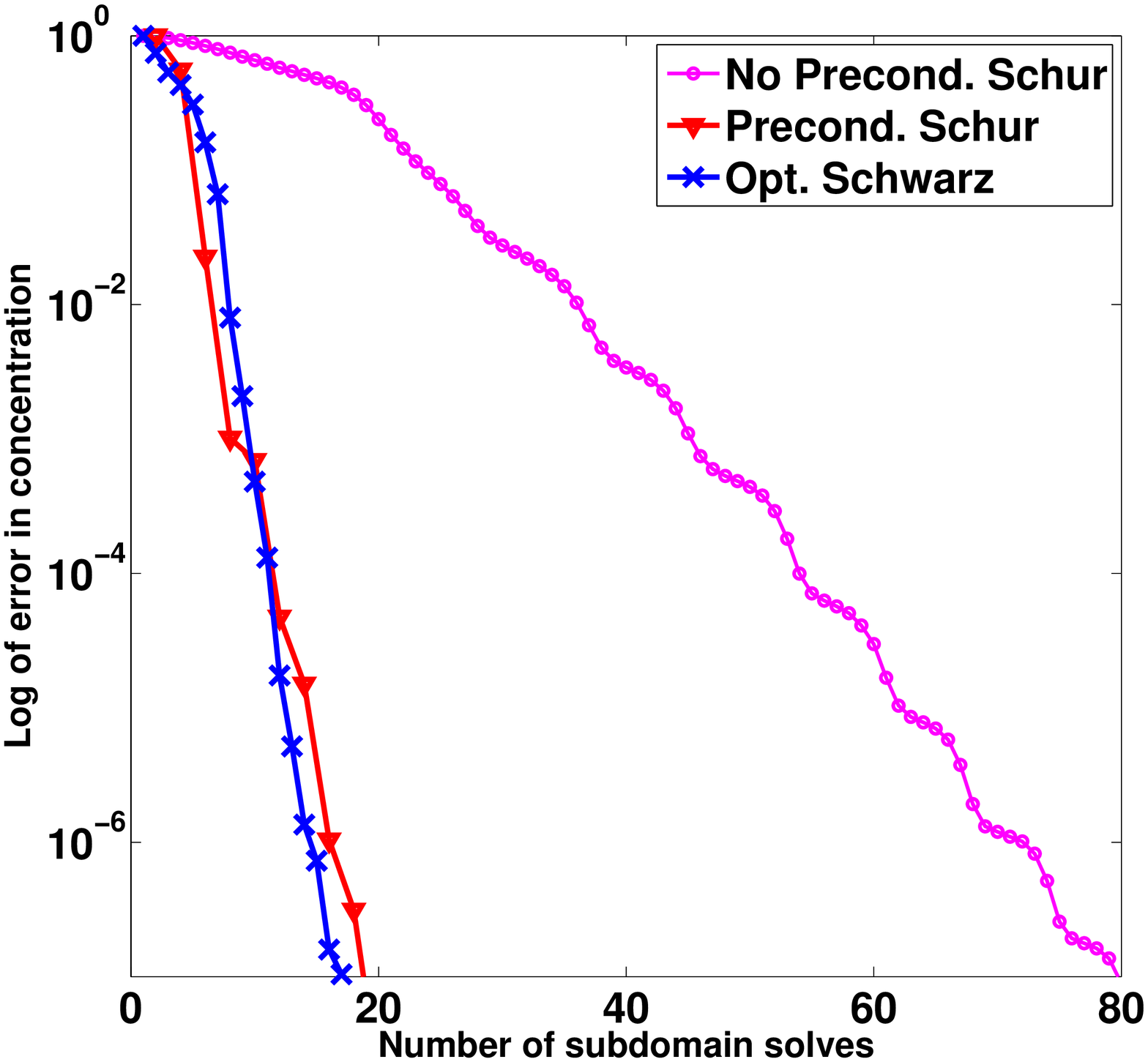} 
\end{minipage} \hspace{16pt}
\begin{minipage}{0.45 \linewidth}
\includegraphics[scale=0.24]{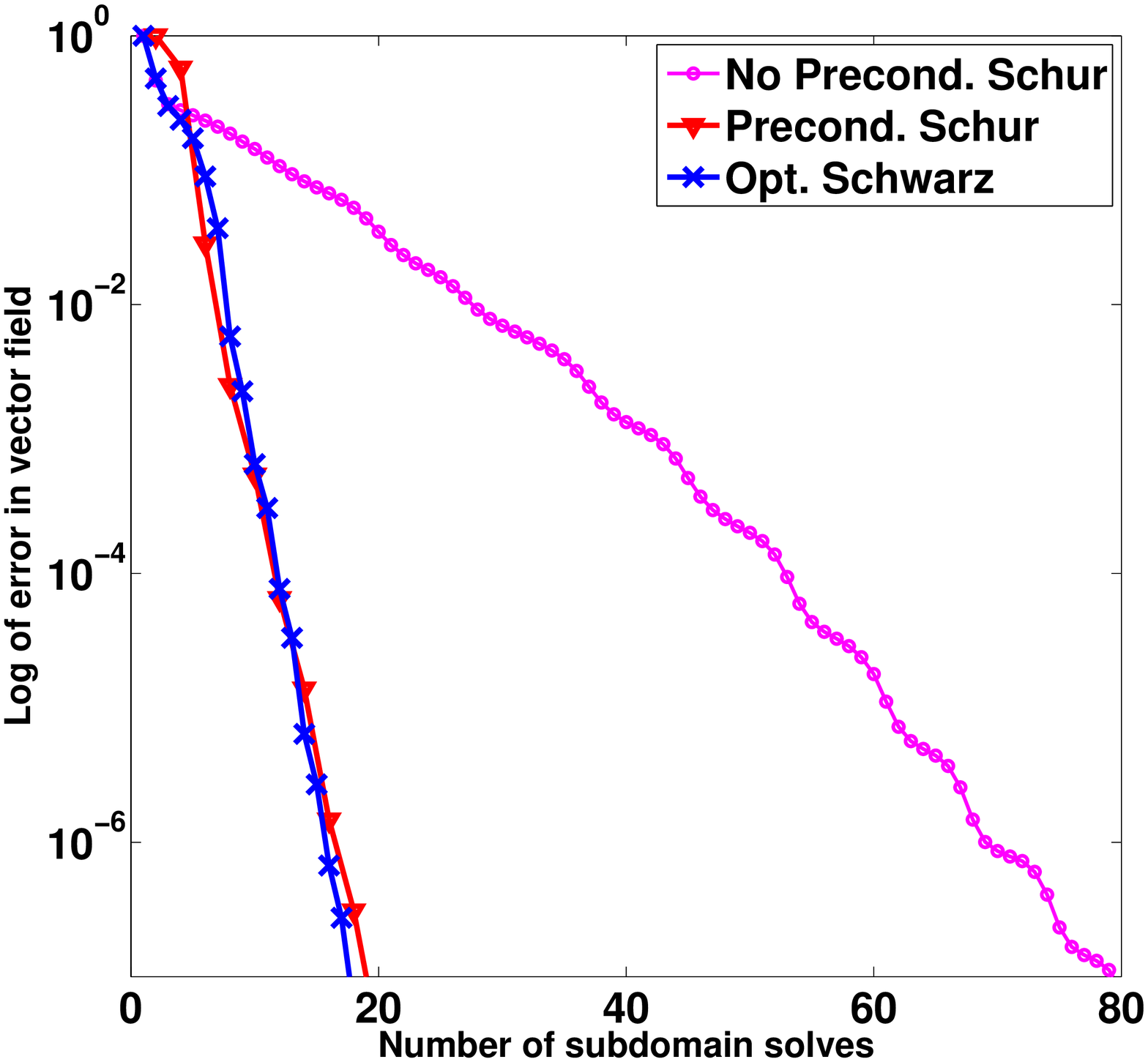} 
\end{minipage} 
\caption{Convergence curves using GMRES: errors in c (on the left) and error in $ \br $
  (on the right).} 	
\label{A2Fig:Test3CEAConvZero}  \vspace{-0.2cm}
\end{figure}

To analyze the convergence behavior of each method, as in the previous test cases, we solve a problem with~$ c_{0} = 0$ (thus $c=0$ and $\br=0$).
We start with a random initial guess on the space-time interface and stop the iteration when the errors both in the concentration $ c $ and in the vector field $ \br $ are less than $ 10^{-6} $ (Figure \ref{A2Fig:Test3CEAConvZero}). We see that the GTP-Schur method with the preconditioner significantly improves the convergence speed compared to the case with no preconditioner, which makes it and the GTO-Schwarz method comparable. This is because the diffusion is dominant in this case. The errors in $ c $ and $ \br $ behave quite similarly. 
%
%

\begin{figure}[H]
\centering
\begin{minipage}{0.45 \linewidth}
\hspace{-1.5cm}
  \includegraphics[scale=0.26]{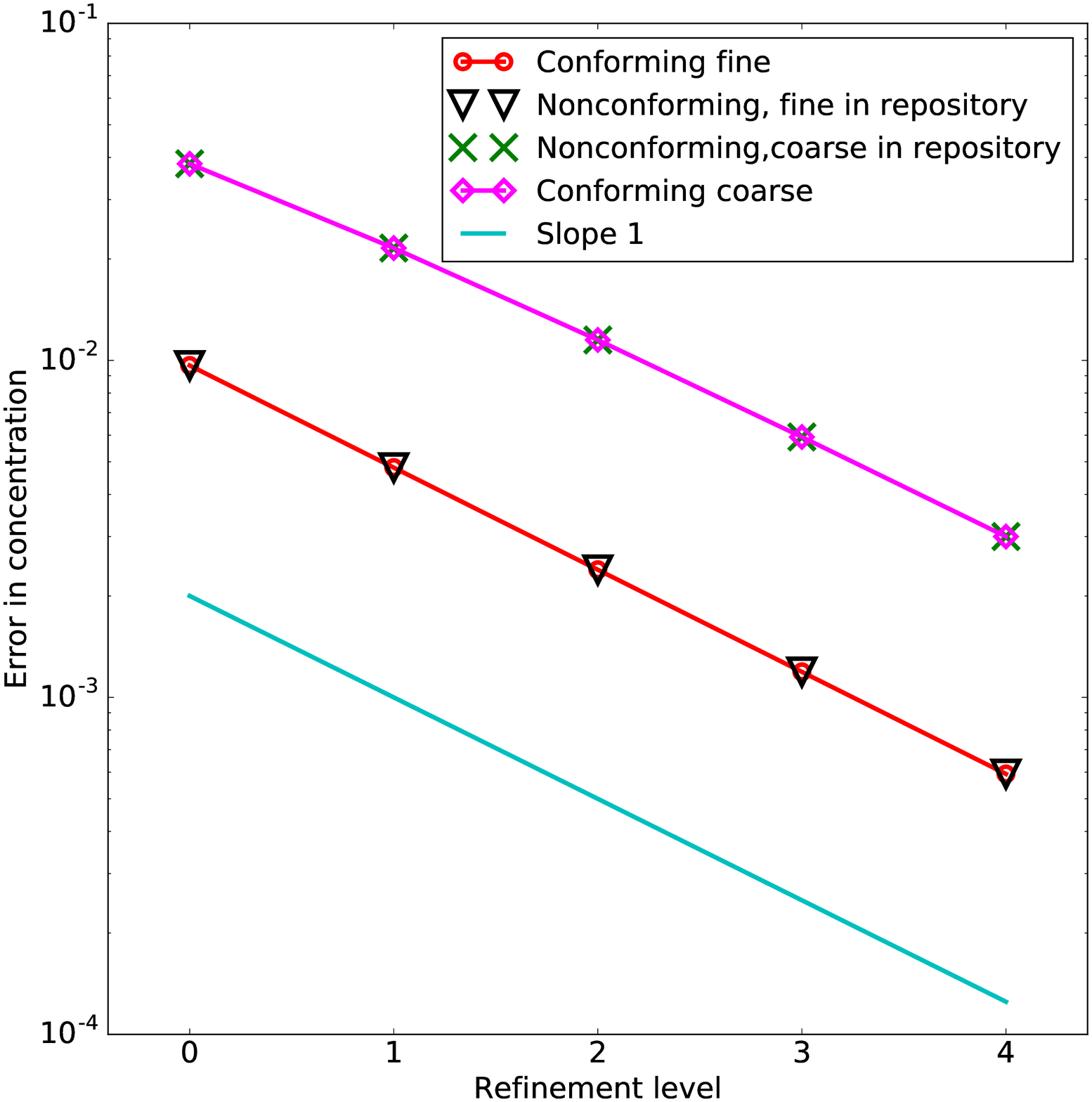} 
\end{minipage} \hspace{10pt}
\begin{minipage}{0.45 \linewidth}
\hspace{-0.5cm}
  \includegraphics[scale=0.26]{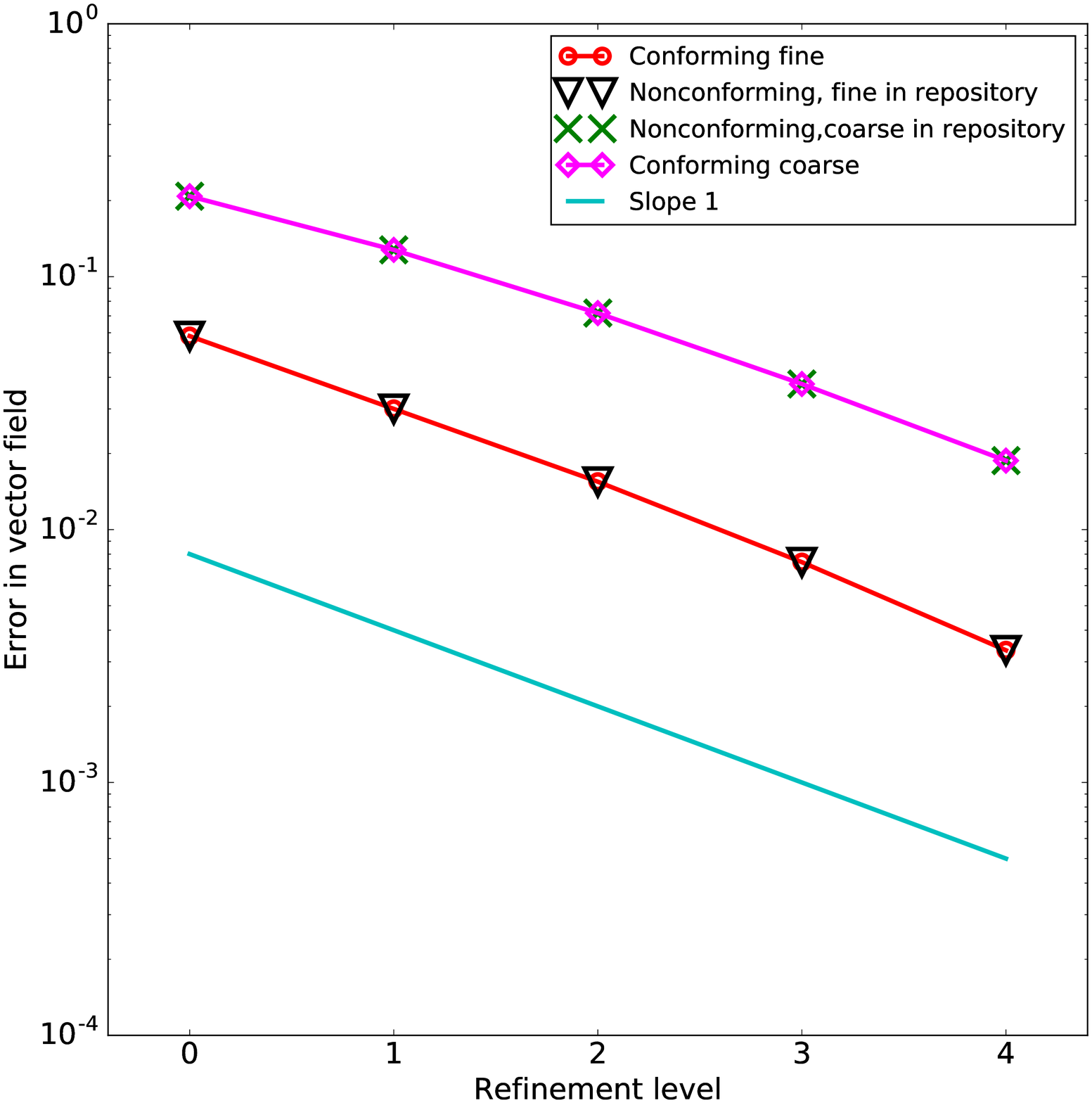} 
\end{minipage} 
\vspace{-0.2cm}
\caption{Errors in $ c $ (left) and $ \br $ (right) in logarithmic scales between the reference and the multidomain solutions versus the refinement level.} 	
\label{A2Fig:Test3CEAConvTime} \vspace{-0.2cm}
\end{figure}

Consider now the initial condition $ c_{0} \neq 0 $ defined in \eqref{A2CEA:IC}. We check to see whether the nonconforming time grids preserve the accuracy in time.  We consider four initial time grids, which we then refine 4 times by a factor of 2, 
\\
$\bullet$ Time grid 1 (conforming fine): $ \Dt_{i} = T/250$, $\forall i$.\\
$\bullet$ Time grid 2 (nonconforming, fine in the repository): $ \Dt_{5} = T/250 $ and
              $ \Dt_{i}~=~T/50$, $\forall i \neq 5 $.\\
$\bullet$ Time grid 3 (nonconforming, coarse in the repository): $ \Dt_{5} = T/50 $ and
              $ \Dt_{i}~=~T/250$,$\forall i \neq 5 $.\\
$\bullet$ Time grid 4 (conforming coarse): $ \Dt_{i} = T/50$, $\forall i$.\\
%
Note that the advection time steps are equal to the diffusion time steps. The time steps are then refined several times by a factor of 2. In space, we fix a conforming rectangular mesh
and we compute a reference solution on a very fine time grid, with $ \Dt = \Dt_{a} =  T/(250 \times 2^{6}) $. Figure~\ref{A2Fig:Test3CEAConvTime} shows the error in the $ L^{2}(0,T; L^{2}(\Omega)) $-norm of the concentration $ c $ and of the vector field $ \br $ versus the refinement level.
%
We observe that first order convergence is preserved in the nonconforming case and the errors obtained in the nonconforming case with a fine time step in the repository (Time grid $ 2 $, in blue)
are nearly the same as in the finer conforming case (Time grid $ 1 $, in red). Thus the use of nonconforming grids (where the ratio of the fine time step to the coarse time step is $ 5 $) preserves the accuracy in time of the monodomain scheme.
%
%
%

\begin{figure}[htbp]
\centering
\begin{minipage}[c]{0.35 \linewidth}
\includegraphics[scale=0.2]{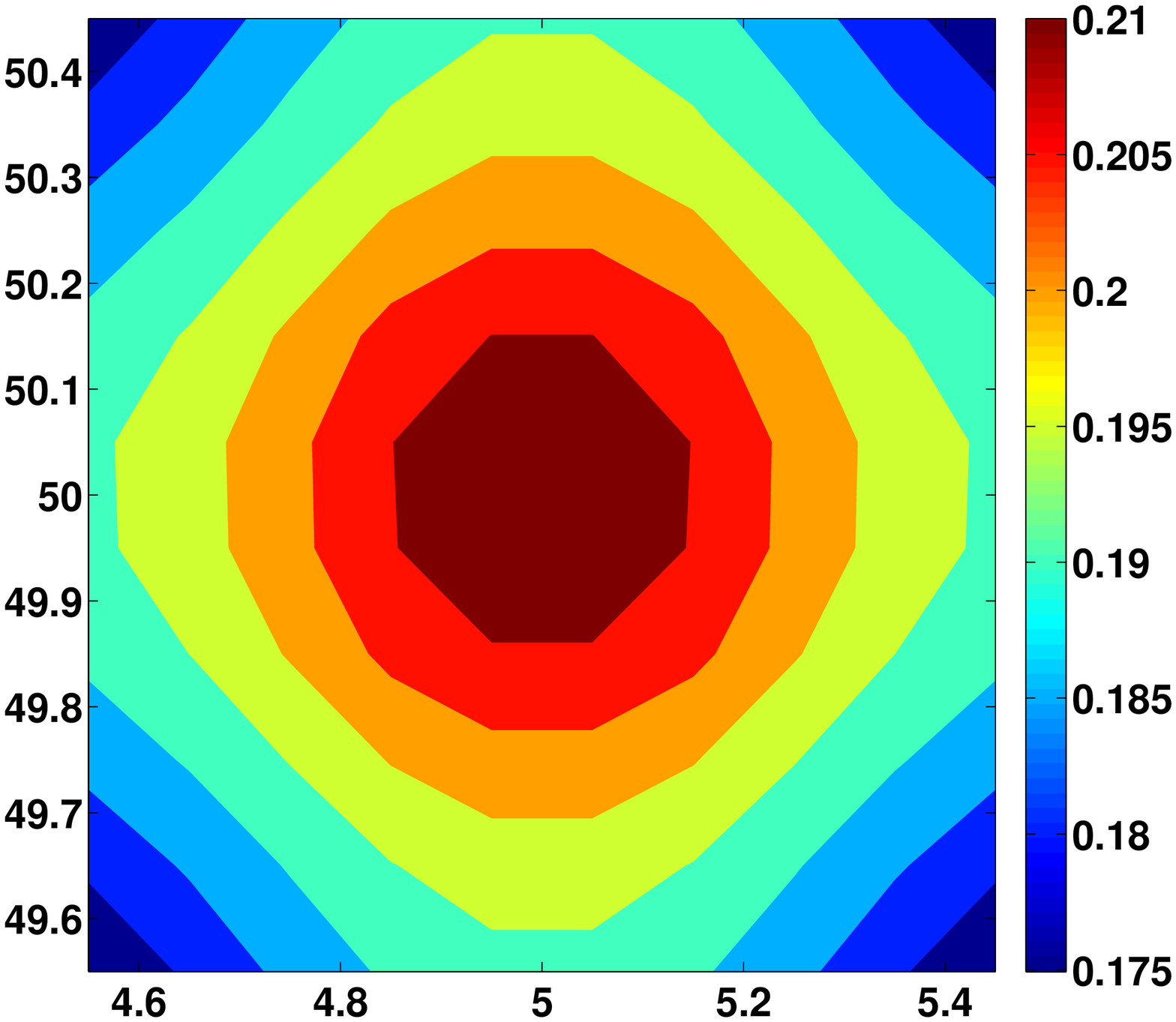} 
\end{minipage} \hspace{-20pt}
\begin{minipage}[c]{0.35 \linewidth}
\includegraphics[scale=0.2]{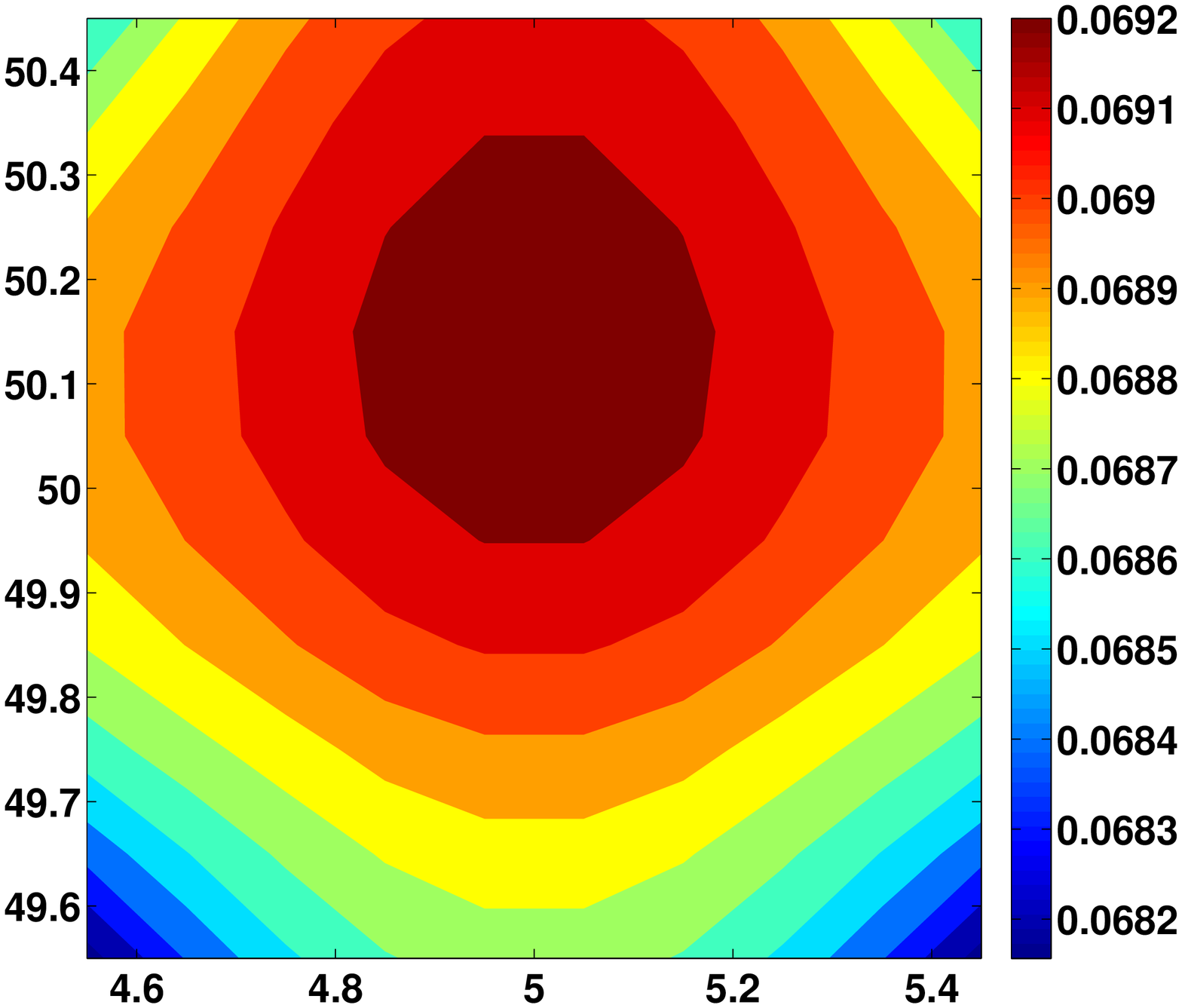}
\end{minipage} \hspace{-20pt}
\begin{minipage}[c]{0.35 \linewidth}
\includegraphics[scale=0.2]{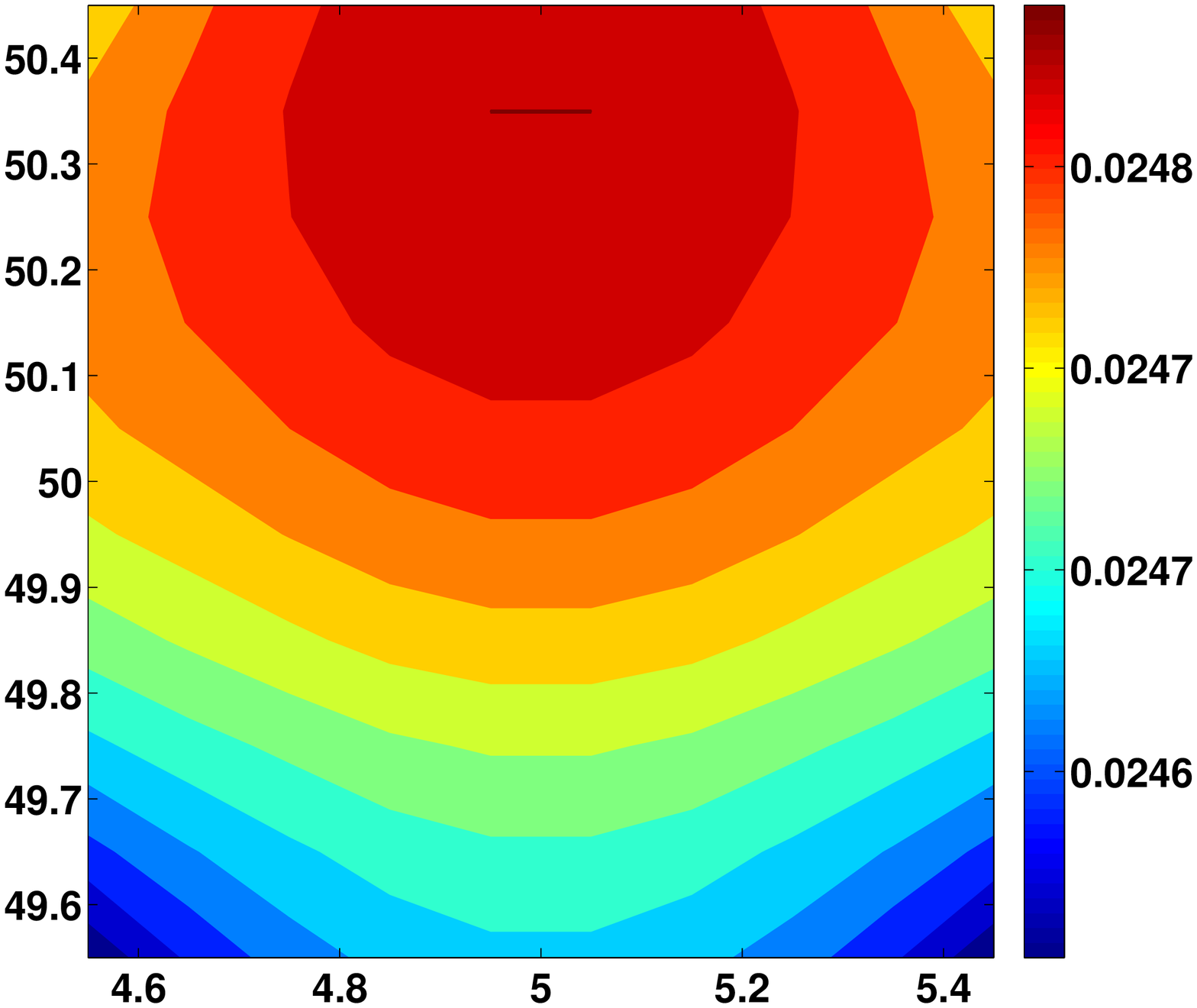}
\end{minipage} 
\\
\begin{minipage}[c]{0.35 \linewidth}
  \includegraphics[scale=0.2]{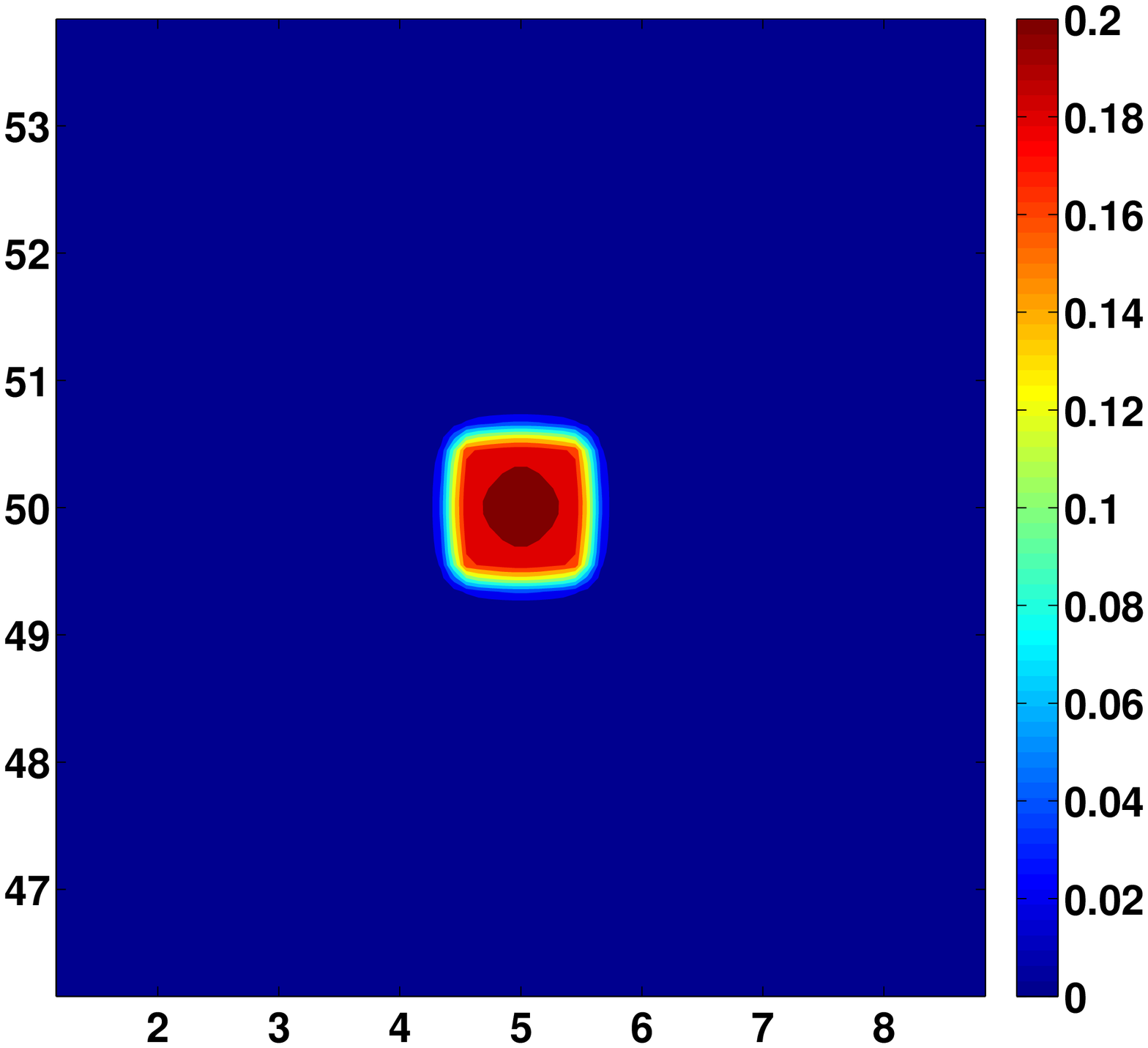}
\end{minipage} \hspace{-20pt}
\begin{minipage}[c]{0.35 \linewidth}
\includegraphics[scale=0.2]{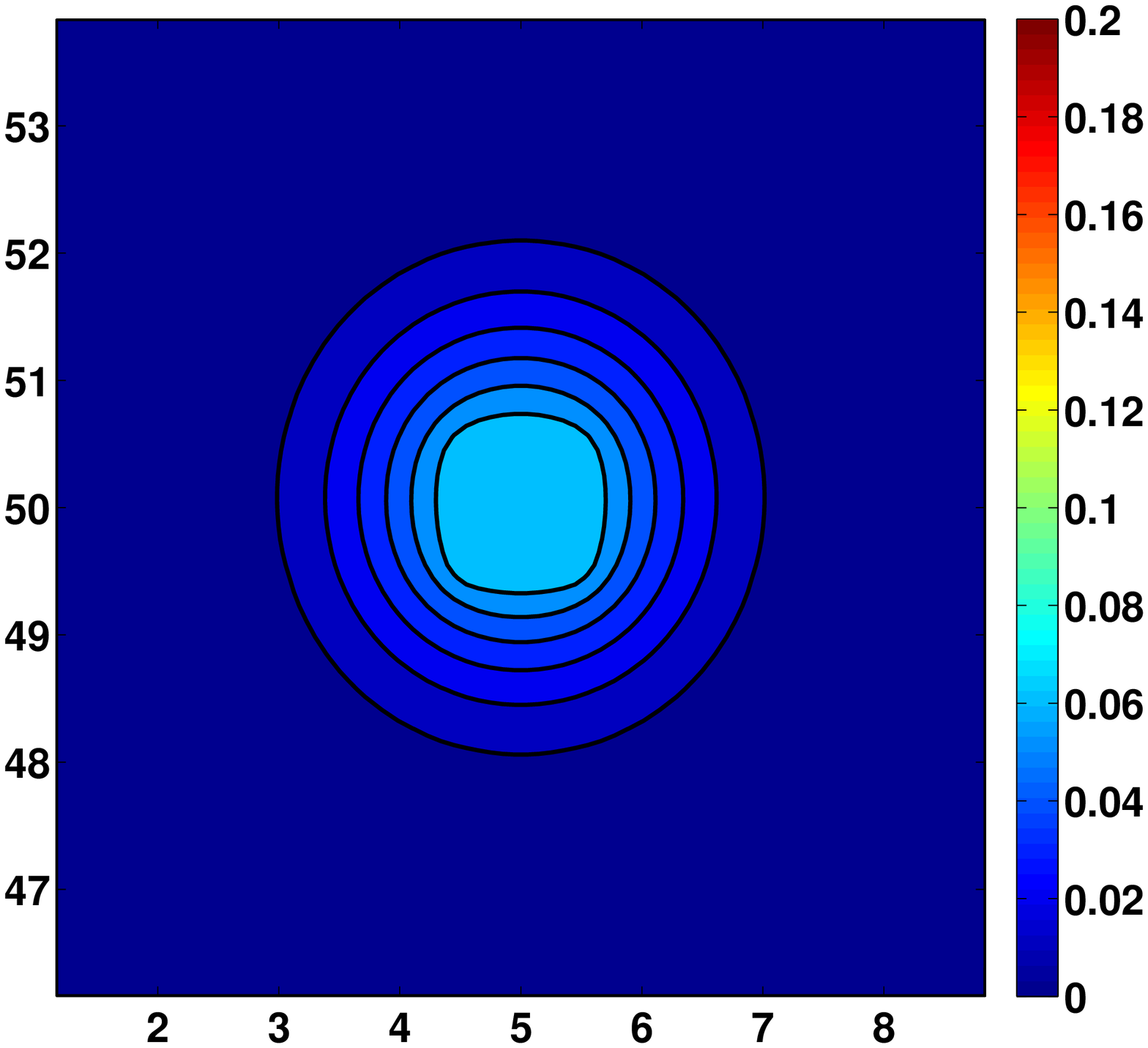}
\end{minipage} \hspace{-20pt}
%
%
%
%
\begin{minipage}[c]{0.35 \linewidth}
\includegraphics[scale=0.2]{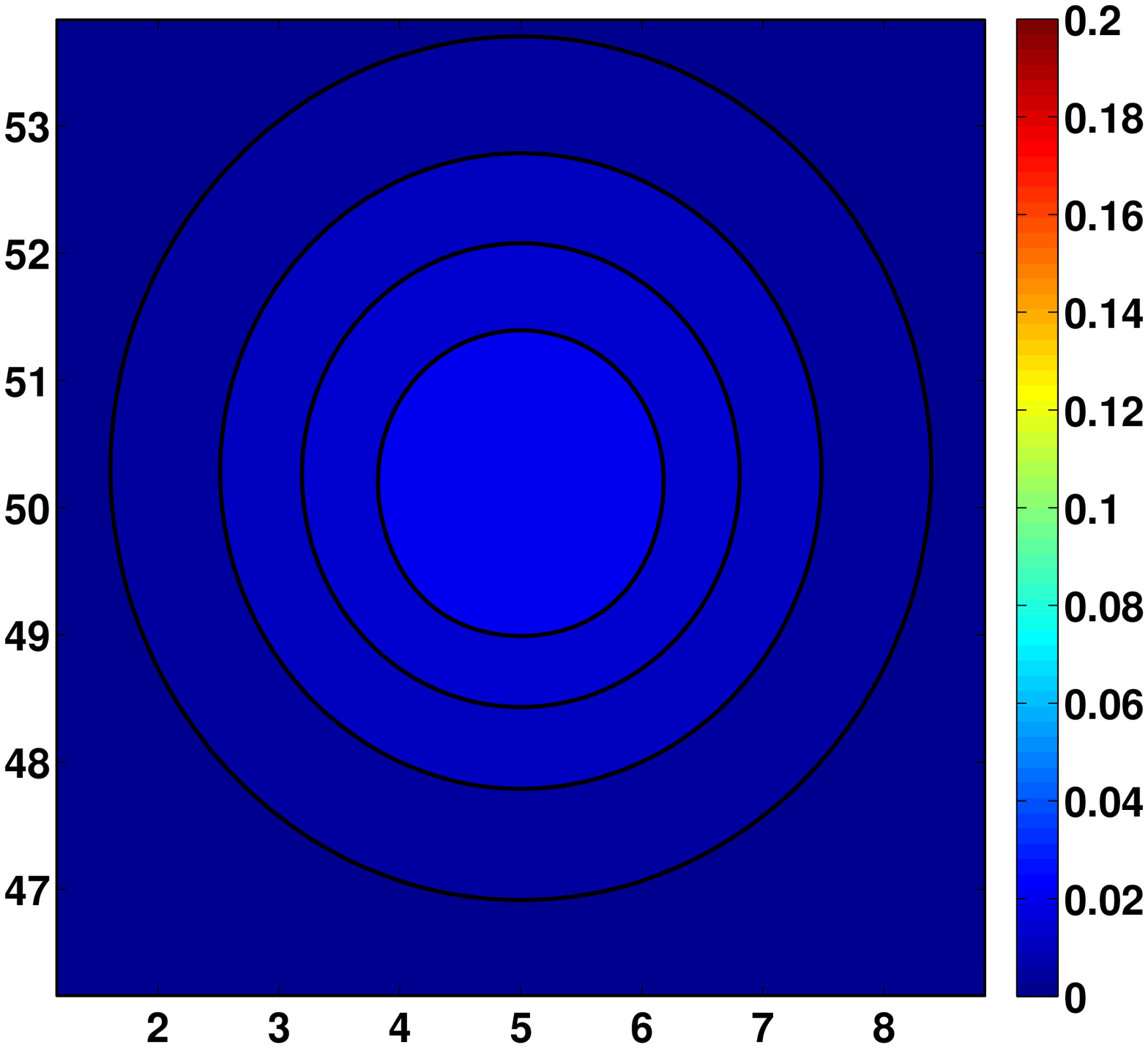}
\end{minipage}
\caption{Snapshots of the concentration in the repository (top) and in the host rock (bottom) after approximately 100 years, 10000 years and 20000 years respectively} 	
\label{A2Fig:CEASolutionTW} \vspace{-0.2cm}
\end{figure}

We now use time windows where the initial guess of the $ (N+1)^{\text{st}} $ time window is calculated from the information at the final time of the $ N^{\text{th}} $ time window (see~\cite{PhuongThesis}), which helps reduce considerably the number of iterations required to reach the same tolerance compared with an arbitrary initial guess. Since the size of the time windows is uniform, we can use the same optimized parameters for all time windows for the GTO-Schwarz method. In each time window, we stop the iterations when the relative residual is less than $ 10^{-3} $. From the observation above, the maximum number of iterations in each time window is not greater than $ 5 $ (equivalent to $ 10 $ subdomain solves) for the GTP-Schur method (with the Neumann-Neumann preconditioner) and  is not greater than $ 5 $ (equivalent to $ 5 $ subdomain solves) for the GTO-Schwarz method.
Figure~\ref{A2Fig:CEASolutionTW} show the concentration in the repository (left) and in the host rock (right) after 1 ($ \approx 100 $ years), 50  ($ \approx 5000 $ years), 100  ($ \approx 10000 $ years) and 200  ($ \approx 20000 $ years) time windows respectively. We use different color scales for the solution in the repository to see clearly the effect of the advection field, while we use same color scales for the solution in the host rock to see the spreading of the contaminant in time. The concentration field behaves as expected and the migration of the radionuclide from the repository to the surrounding medium takes place very slowly.
%
%
%
%
%
	\subsection{A simulation for a surface, nuclear waste storage}
	\label{A2subsec:ANDRA}
\begin{figure}[h]
\centering
\begin{minipage}[c]{0.9 \linewidth}
\begin{center}
\includegraphics[scale=0.4]{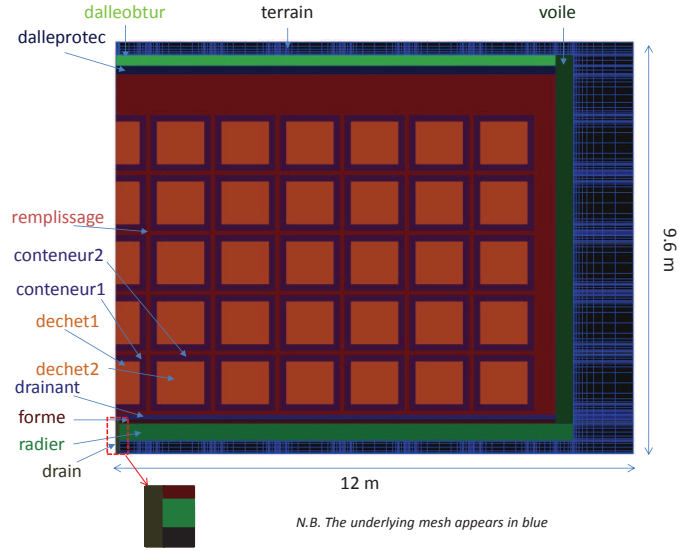} 
\end{center}
\end{minipage}  \vspace{-0.2cm}
\caption{The geometry of the test case.} 	
\label{A2Fig:ANDRAdomain} \vspace{-0.2cm}
\end{figure}

%
%
We consider a test case designed by ANDRA for a surface storage of nuclear waste of short half-life. The computational domain is depicted in Figure~\ref{A2Fig:ANDRAdomain} with different physical zones, where the waste is stored in square boxes (\textit{d\'echet} zone). The properties of these zones are given in Table~\ref{A2Tab:ANDRAdata}. Note that in our calculation, we use the effective diffusion, defined by $ d_{\text{eff}} = \phi \times d_{\text{m}} $. The advection field is governed by Darcy's law together with the law of mass conservation~\eqref{A2CEADarcyflow}, where we have used
the hydraulic head field $h$ instead of the pressure field $p$.
Dirichlet conditions are imposed on top, $h= 10$m and on bottom $ h = 9.998$m of the domain and no-flow boundary on the left and right sides.
The final time is $ T_{f} = 500 $ years. The source term is $ f = 0 $ and the initial condition is such that 
\begin{equation*} \label{A2ANDRA:IC}
c_{0}= \left \{ \begin{array}{ll} 1, & \text{in \textit{d\'echet1} and \textit{d\'echet2}}, \\
0, & \text{elsewhere}.
\end{array} \right .
\end{equation*}
Boundary conditions of the transport problem are homogeneous Dirichlet conditions on top and bottom, and homogeneous Neumann conditions on the left and right hand sides. 

\begin{figure}[H]
\vspace{0.2cm}
\begin{minipage}{0.4 \linewidth}
\hspace{-1cm}\includegraphics[scale=0.26]{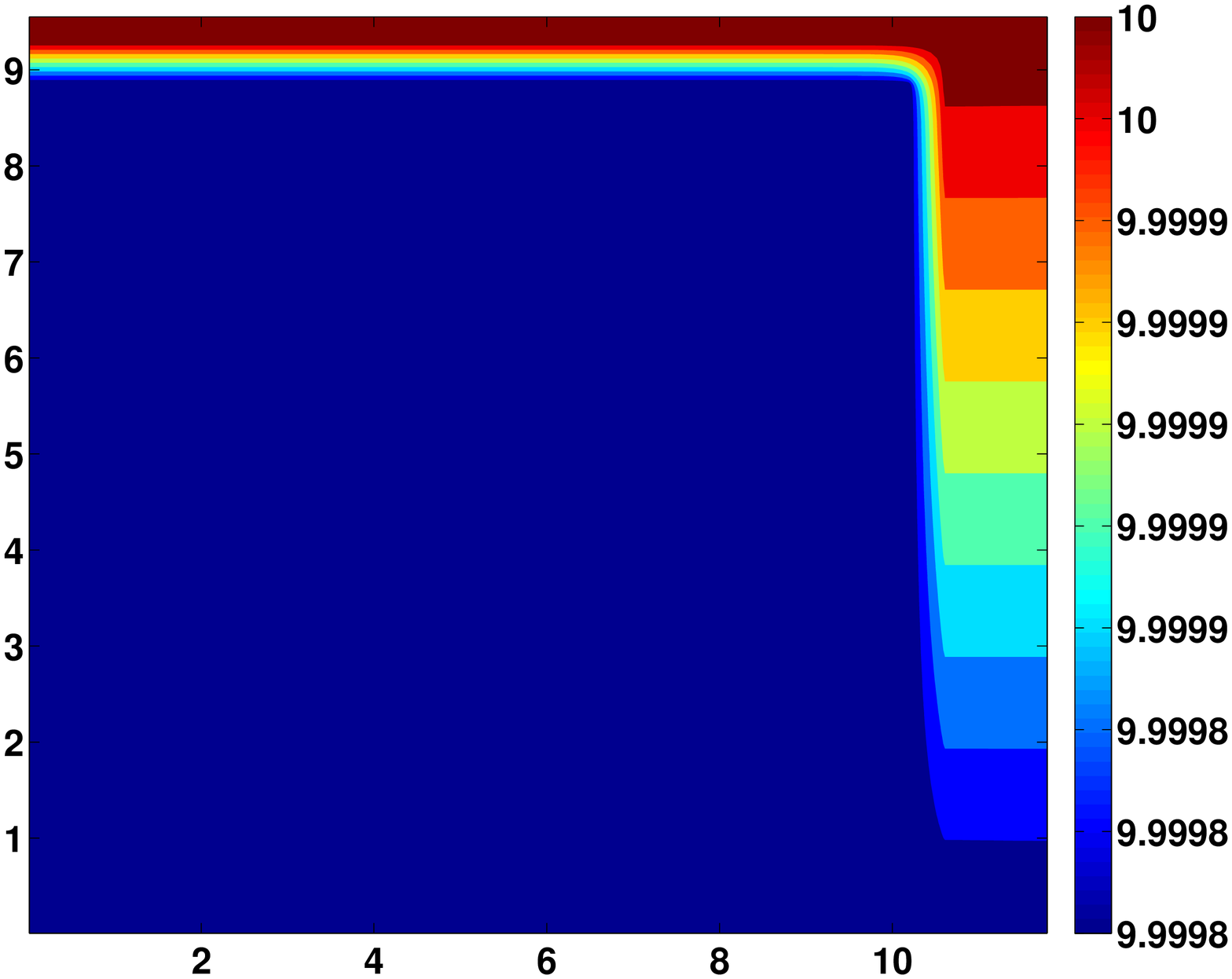} 
\end{minipage} \hspace{4pt}
\begin{minipage}{0.5 \linewidth}
\vspace{-0.1cm}
\setlength{\unitlength}{1pt} 
\begin{picture}(140,100)(0,0)
\thicklines
\put(45,-40){\includegraphics[height=6.cm]{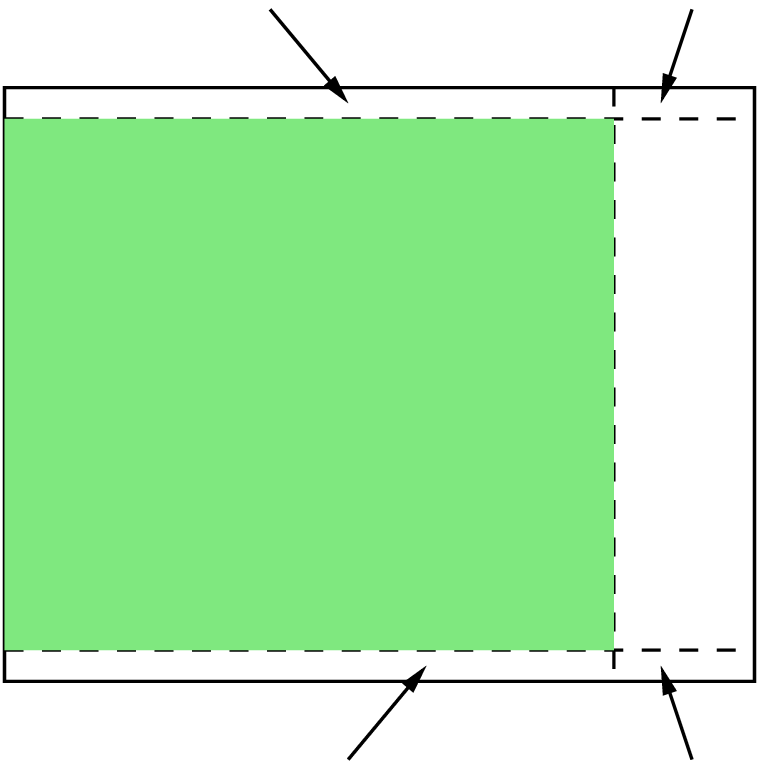} \\}
\put(107,-38){$ \Omega_{1} $}
\put(205,-38){$ \Omega_{2} $}
\put(120,45){$ \Omega_{3}$}
\put(195,45){$ \Omega_{4} $}
\put(100,135){$ \Omega_{5} $}
\put(200,135){$ \Omega_{6}$}
\end{picture}
\end{minipage} 
\caption{The hydraulic head field and the decomposition of the domain.} 	
\label{A2Fig:ANDRADarcyDecom}   \vspace{-0.2cm}
\end{figure}

For the spatial discretization (for both the flow and transport equations), we use a non-uniform rectangular mesh, shown in Figure~\ref{A2Fig:ANDRAdomain} in blue, with $ 171 $ cells in the $ x-$direction and $ 158 $ cells in the $ y- $direction. The mesh size is $ \Delta x \approx 0.42$m.  The hydraulic head is approximated by mixed finite elements and is shown in Figure~\ref{A2Fig:ANDRADarcyDecom} (left). We decompose the domain into $ 6 $ rectangular subdomains in such a way that the black zone (\textit{terrain}) is separated from the rest and subdomain $ \Omega_{3} $ includes the \textit{dallerobtur}, \textit{voile}, \textit{radier} and a part of \textit{drain} zones (see Figure~\ref{A2Fig:ANDRADarcyDecom} (right)).
The transport is dominated by diffusion in subdomain $ \Omega_{3} $ (the maximum of the local P\'eclet number $\text{Pe}_{L} \approx  0.0032 $) and is dominated by advection (with $\text{Pe}_{L} \approx  2.75 $) in the other subdomains. The time steps due to the CFL condition are $ \Dt_{a,3} \leq 0.6551 $ years and very small elsewhere, $ \Dt_{a,i} \leq 6.0874 \, 10^{-5} $ years, $\forall i \neq 3 $.
\begin{table}[h]
\centering
{\footnotesize
\begin{tabular}{|l|l|l|l|l|}
  \hline
	    		Zone               						& Hydraulic conductivity       & Porosity                    & Molecular diffusion  				 \\ 
	    													& $ K $ (m/year) 				  & $ \phi $					& $ d_{\text{m}}$ (m$^{2}$/year)   \\ \hline
	    		terrain								    & $ 94608 $   									& $ 0.30 $                  & $ 1 $    	 					\\ \hline 
	    		apron (\textit{radier})									& $ 3.1536 \, 10^{-4} $   					& $ 0.15 $                  & $ 6.31 \, 10^{-5}$     	  \\ \hline 	
	    		shape (\textit{forme})									& $ 3.1536 \, 10^{-3} $   					& $ 0.20 $                  & $ 1.58 \, 10^{-3} $     	  \\ \hline 
	    		draining area (\textit{drainant})								& $ 94608 $   									& $ 0.30 $                  & $ 5.36 \, 10^{-2} $     	  \\ \hline
	    		veil (\textit{voile})									& $ 3.1536 \, 10^{-3} $   					& $ 0.20 $                  & $ 1.58 \, 10^{-3} $   	  \\ \hline
	    		filling (\textit{remplissage})							& $ 5045.76 $   								& $ 0.30 $                  & $ 5.36 \, 10^{-2} $     	  \\ \hline
	    		protection slab (\textit{dalleprotec})							& $ 3.1536 \, 10^{-3} $   					& $ 0.20 $                  & $ 1.58 \, 10^{-3} $     	  \\ \hline
	    		closing slab (\textit{dallerobtur})							& $ 3.1536 \, 10^{-3} $ 					& $ 0.20 $                  & $ 1.58 \, 10^{-3} $     	  \\ \hline
	    		drain 									& $ 94608 $   									& $ 0.30 $                  & $ 1 $     	  					\\ \hline
	    		container1/container2 (\textit{conteneur1/conteneur 2})		& $ 3.1536 \, 10^{-4} $   					& $ 0.12 $                  & $ 4.47 \, 10^{-4} $     	  \\ \hline
	    		waste1/waste2 (\textit{d\'echet1/d\'echet2})				& $ 3.1536 \, 10^{-4} $   					& $ 0.30 $                  & $ 1.37 \, 10^{-3} $     	  \\ \hline	     		
\end{tabular}}
\caption{Data for flow and transport problems.}
\label{A2Tab:ANDRAdata}  \vspace{-0.2cm}
\end{table}
As in Section~\ref{A2subsec:CEA}, we use time windows with size $ T = 10 $ years. We consider the first time window, $ (0,T) $ and use nonconforming time grids with $ \Dt_{3} =0.1$ years and $ \Dt_{i} = 0.5 $ years, $\forall i \neq 3 $. The advection steps, satisfying the CFL conditions, are $ \Dt_{a,3} = \Dt_{3} $ and $ \Dt_{a,i} = 825 \,\Dt_{i}, \, \forall i \neq 3 $.
%
%
We use a zero initial guess on the space-time interface, and perform GMRES for both the GTP-Schur method (with Neumann-Neumann preconditioner) and the GTO-Schwarz method. We haven't shown the unpreconditioned GTP-Schur as it converged more slowly than the preconditioned version. We compute the errors of the difference between the multidomain solution and a reference solution computed with a very fine, conforming time grid.  
\begin{figure}[H]
\centering
\begin{minipage}{0.5 \linewidth}
\includegraphics[scale=0.26]{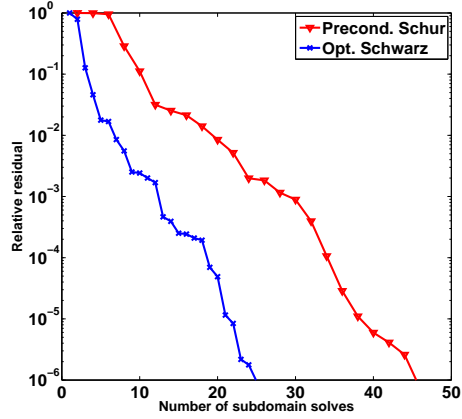} 
\end{minipage} 
\vspace{-0.2cm}
\caption{Relative residuals of GMRES for the GTP-Schur method (with the Neumann-Neumann preconditioner) and the GTO-Schwarz method.} 	
\label{A2Fig:ANDRAStoppingIterRelres}    \vspace{-0.2cm}
\end{figure}
%
%
We show in Figure~\ref{A2Fig:ANDRAStoppingIterRelres} the relative residuals for each method versus the number of subdomain solves.
We observe that the GTO-Schwarz method converges about twice as fast as the GTP-Schur method. This is a case where there is sufficient advection to make GTO-Schwarz faster than GTP-Schur with Neumann-Neumann preconditioning, but the advection term is small enough for the Neumann-Neumann preconditioner to be effective. The errors in $c$ and
$\br$ due to the schemes are obtained for both methods when the relative residual is less than $ 10^{-2} $.

\begin{figure}[H]
\begin{flushleft}
\begin{minipage}{0.5 \linewidth}
\includegraphics[scale=0.22]{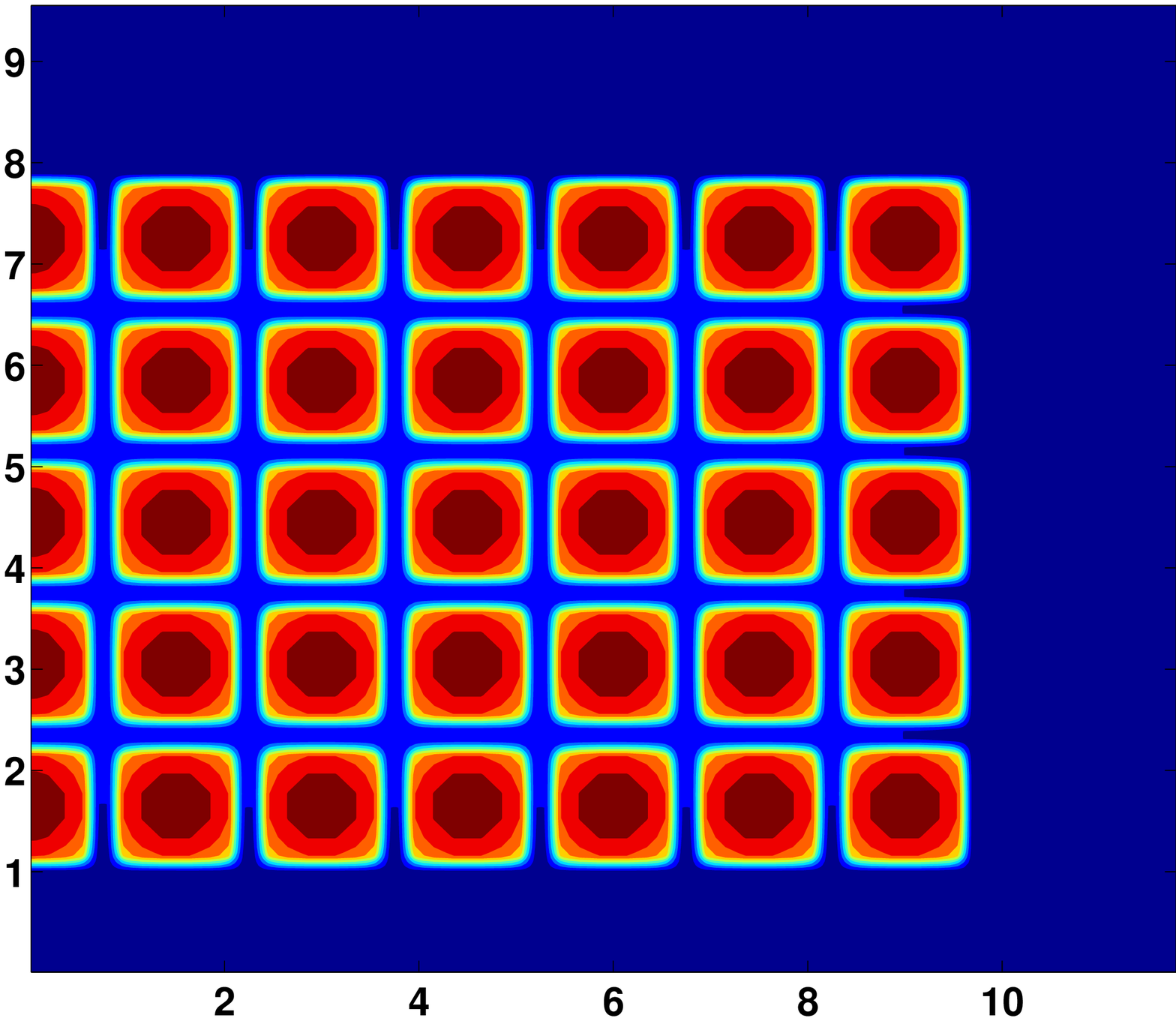} 
\end{minipage} \hspace{2pt}
\begin{minipage}{0.4 \linewidth}
\includegraphics[scale=0.22]{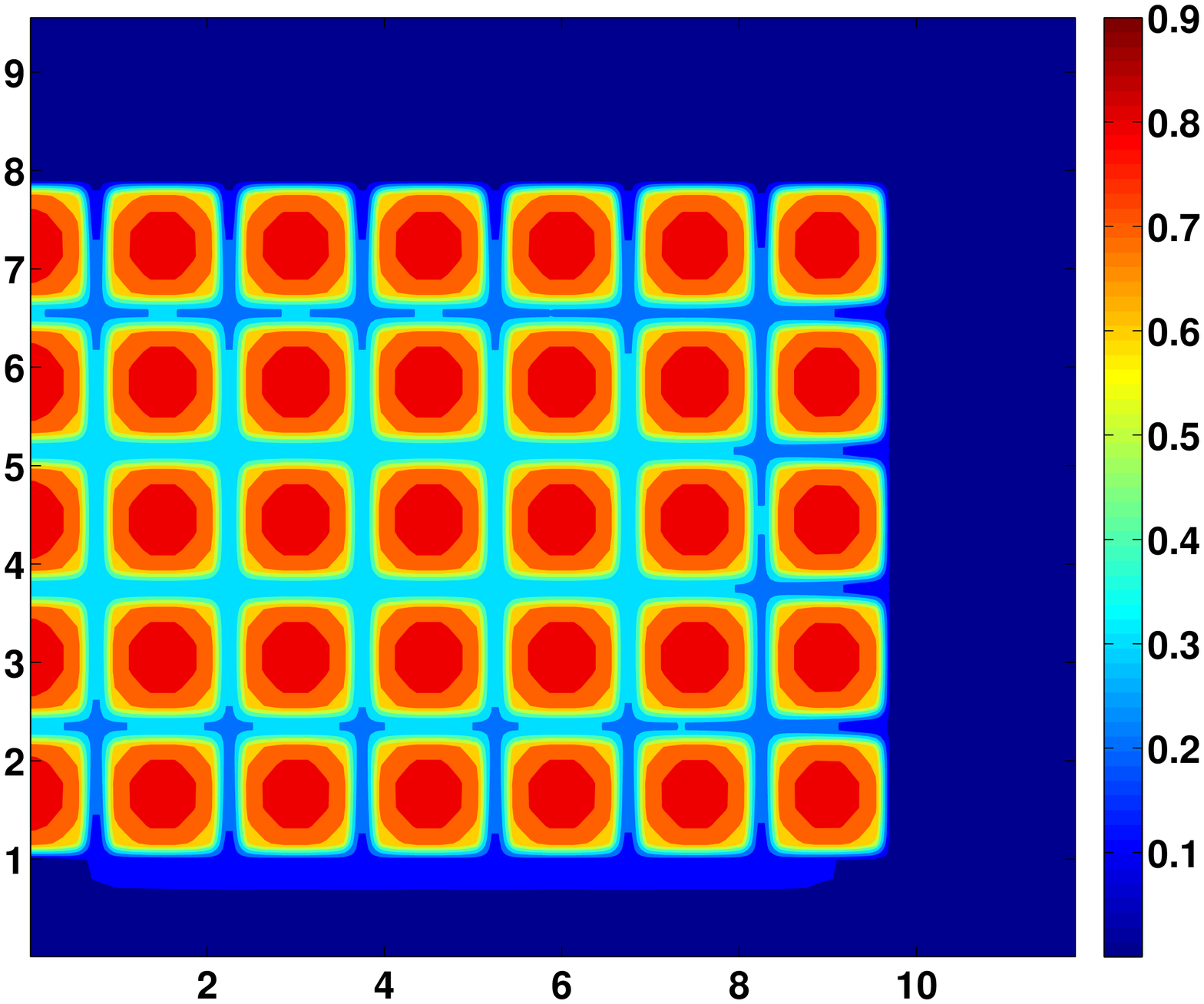} 
\end{minipage} 
\begin{minipage}{0.5 \linewidth}
  \vspace{-0.6cm}
\includegraphics[scale=0.22]{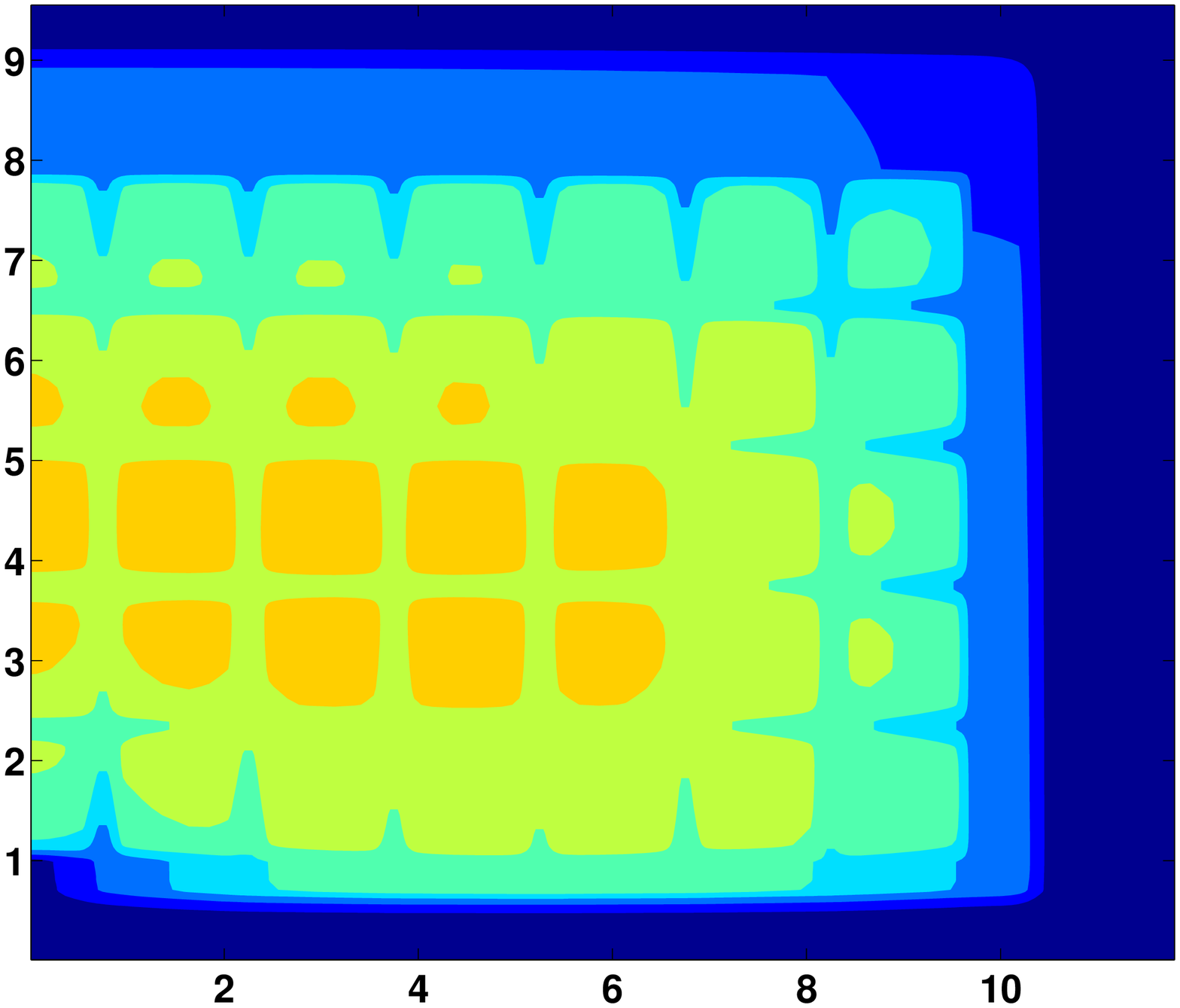} 
\end{minipage} \hspace{2pt}
\begin{minipage}{0.4 \linewidth}
  \vspace{-0.6cm}
\includegraphics[scale=0.22]{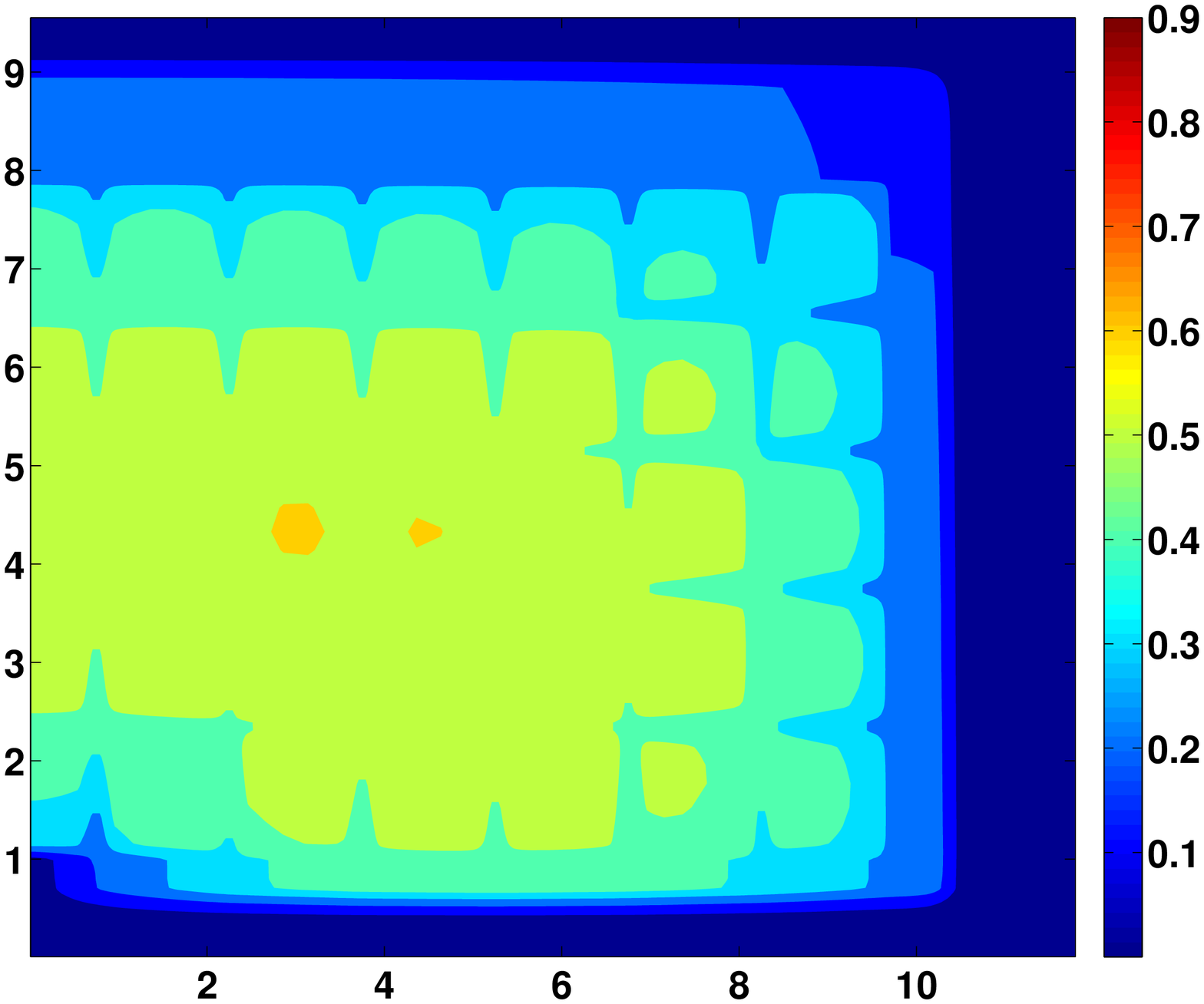} 
\end{minipage} 
\end{flushleft}
\vspace{-0.2cm}
\caption{Snapshots of the concentration after 20~years, 50~years, 350~years and 500~years respectively.} 	
\label{A2Fig:ANDRASolutionTW} \vspace{-0.2cm}
\end{figure}

Next, we run the two methods for $ 50 $ time windows and stop the iterations in each time window when the relative residual is less than $ 10^{-2} $. The maximum number of iterations in each time window is not greater than $ 10 $ (equivalent to $ 20 $ subdomain solves) for the GTP-Schur method (with the Neumann-Neumann preconditioner) and  is not greater than $ 8 $ (equivalent to $ 8 $ subdomain solves) for the GTO-Schwarz method. 
Figure~\ref{A2Fig:ANDRASolutionTW} shows the concentration field after 2 ($20$~years), 5  ($ 50 $~years), 35  ($ 350 $~years) and 50  ($ 500 $~years) time windows respectively. We see that the radionuclide escapes from the waste packages and slowly migrates into the surrounding area. Due to the specific design of the storage and under the effect of advection, the radionuclide tends to move toward the bottom right corner.
\begin{remark}
In this paper we have not studied parallel implementation or scalability issues. A theoretical study along the lines of~\cite{Gropp:1988:CPI}, and \cite{OSWRDG2}  with nonconforming time grids, could be carried out. Without going into details, we note that the scalability properties of the methods would be very similar to those of a domain decomposition method in space only, with one important difference being that for a space-time method, communication between subdomains only occurs at the end of the time interval. Note also that for a large number of subdomains, one should use a coarse grid or a coarse space correction to preserve numerical scalability.
\end{remark}  
\section*{Conclusion}

In the context of operator splitting, we have extended the two methods derived in~\cite{PhuongSINUM} for the pure diffusion problem to the heterogeneous advection-diffusion problem. Two discrete interface problems corresponding to the generalized time-dependent Steklov-Poincar\'e operator and the OSWR approach with operator splitting have been formulated in a way such that they are equivalent to the discrete monodomain problem and that they enable different advection and diffusion time steps in the subdomains. For the GTP-Schur method, a generalized Neumann-Neumann preconditioner is considered and is validated for different test cases in 2D experiments. Numerical results show that the GTO-Schwarz method outperforms the GTP-Schur method (with or without preconditioner) in terms of subdomain solves needed to reach a fixed error reduction in the solution (by a factor of $ 2 $ to $ 2.5 $ in our test cases). Due to the use of the optimized Robin parameters, the GTO-Schwarz method is robust in the sense that it handles well and consistently both the advection-dominated and diffusion-dominated problems. The GTP-Schur method with the Neumann-Neumann preconditioner works well and converges faster than without a preconditioner when the diffusion is dominant and it also efficiently deals with the case with large jumps in the diffusion coefficients. When the advection is dominant, the Neumann-Neumann preconditioner converges slower than when there is no preconditioner. However, asymptotically the convergence of the Neumann-Neumann preconditioned the GTP-Schur method has a weak dependence on the mesh size of the discretizations while that of the GTP-Schur method with no preconditioner significantly depends on the mesh size. For the GTO-Schwarz method, because of the optimized parameters, which play in some sense the role of a preconditioner, the convergence is weakly dependent on the discretization parameters. In addition, both methods preserve the accuracy in time when nonconforming time steps are used, both for two subdomains and for multiple subdomains: the error due to the nonconforming time grid (with finer time steps in the zones where the solution varies most, i.e. with larger advection and diffusion coefficients) is close to that of the conforming fine grid.
Note that it is known that in the stationnary case the Neumann-Neumann preconditioner performs significantly worse in the presence of advection, whereas the optimization inherent to the optimized Schwarz method makes it insentitive to the advection~\cite{Achdou:DDnonsym:1999,JaphetNatafRogier}.

Two test cases for the simulation of nuclear waste disposal are implemented using nonconforming time grids and time windows. As the geometry of the computational domain is complex and the physical coefficients are highly variable, we consider multiple subdomains.  For this application where the diffusion is dominant, the Neumann-Neumann preconditioned the GTP-Schur method and the GTO-Schwarz method work well but the GTO-Schwarz method converges faster than the GTP-Schur method. We also observe that with an adapted initial guess calculated from the previous time window, one performs only a few iterations in each time window to reach the scheme error. 

\section*{Aknowledgments}
We would like to thank Marc Leconte from ANDRA for his assistance in setting up the
example of Section~\ref{A2subsec:ANDRA}.

We would like to thank the referees for valuable suggestions.
\section*{References}

\bibliography{Advectionpaper_revised}

\end{document}